\documentclass[11pt]{amsart}
\usepackage{amssymb}
\usepackage{epsfig}
\usepackage{graphics}
\usepackage{graphicx}
\renewcommand{\a}{\alpha}

\renewcommand{\d}{\delta}

\newcommand{\e}{\epsilon}

\renewcommand{\l}{\lambda}
\newcommand{\s}{\sigma}

\newcommand{\p}{\phi}

\newcommand{\SJ}{{\mathcal{J}}}

\newcommand{\SL}{{\mathcal{L}}}
\newcommand{\SM}{{\mathcal{M}}}

\newcommand{\SP}{{\mathcal{P}}}

\newcommand{\Z}{\mathbb{Z}}
\newcommand{\C}{\mathbb{C}}

\newcommand{\N}{\mathbb{N}}
\newcommand{\R}{\mathbb{R}}

\newcommand{\CP}{\mathbb{CP}}
\newcommand{\RP}{\mathbb{RP}}

\newcommand{\iso}{\cong}

\newcommand{\id}{\textbf{\textit{I}}}

\newcommand{\Id}{\operatorname{Id}}

\newcommand{\sref}[1]{\sigma^{\text{\rm \scriptsize ref}}_{k,#1}}
\newcommand{\bd}{\partial}
\newcommand{\bbd}{\bar{\partial}}

\newtheorem{proposition}{Proposition}[section]
\newtheorem{theorem}[proposition]{Theorem}
\newtheorem{definition}[proposition]{Definition}
\newtheorem{lemma}[proposition]{Lemma}
\newtheorem{conjecture}[proposition]{Conjecture}

\newtheorem{remark}[proposition]{Remark}

\title{Lagrangian submanifolds and Lefschetz pencils.}

\subjclass[2000]{Primary: 53D12. Secondary: 53D35.}
\date{June, 2004.}
\keywords{Symplectic, Lefschetz pencil, Lagrangian submanifold}

\author{Denis Auroux}
\address{Department of Mathematics \\ M.I.T.\\
Cambridge, MA 02139-4307 \\ U.S.A.}
\email{auroux@math.mit.edu}

\author{Vicente Mu\~noz}
\address{Instituto de Matem\'aticas y F\'{\i}sica Fundamental\\
Consejo Superior de Investigaciones Cien\-t\'{\i}\-fi\-cas \\ C/ Serrano 113-bis \\ 28006 Madrid \\
Spain} \email{vicente.munoz@imaff.cfmac.csic.es}

\author{Francisco Presas}
\address{Departamento de Matem\'aticas \\
Universidad Aut\'onoma de Madrid
\\ 28049 Madrid \\ Spain}
\email{francisco.presas@uam.es}

\begin{document}

\renewcommand{\theenumi}{\roman{enumi}}

\begin{abstract}
 Given a Lagrangian submanifold in a symplectic manifold and a
 Morse function on the submanifold, we show that there is an
 isotopic Morse function and a symplectic Lefschetz pencil on the
 manifold extending the Morse function to the whole manifold. {}From
 this construction we define a sequence of symplectic invariants
 classifying the isotopy classes of Lagrangian spheres in a
 symplectic $4$-manifold.
\end{abstract}

\maketitle

\section{Introduction} \label{sec:introduction}

For a symplectic manifold $(M,\omega)$, S. Donaldson has proved in
\cite{Do00} the existence of symplectic Lefschetz pencils using
the recently introduced asymptotically holomorphic techniques
\cite{Do96,Au97}.

To give the definition of a Lefschetz pencil, recall that a chart
$(\psi,U)$, $\psi=(z_1,\ldots, z_n): U\subset M \to \C^{n}$ is
adapted at a point $p\in M$ if $\psi(p)=0$ and $\psi^*(J_0)$ is
tamed by $\omega$ (where $J_0$ is the standard complex structure
on $\C^n$); equivalently, this means that complex lines in the
local coordinates are symplectic with respect to $\omega$.

\begin{definition}\label{def:pencil}
A symplectic Lefschetz pencil associated to a symplectic manifold
$(M,\omega)$ consists of the following data:
\begin{enumerate}
 \item A codimension $4$ symplectic submanifold $N$.
 \item A surjective map $\p:M-N \to \CP^1$.
 \item A finite set of points $\Delta\subset M-N$ away
 from which the map $\p$ is a submersion.
\end{enumerate}
Moreover the data satisfy the following local models
\begin{enumerate}
 \item For any point $p\in N$, there exists an adapted chart $(z_1,
 \ldots, z_n)$ for which the submanifold $N$ has local equation
 $\{z_1=z_2=0 \}$ and such that $\p(z_1, \ldots, z_n)=z_2 / z_1$.
 \item For any point $p\in \Delta$, there exists an adapted
 chart $(z_1, \ldots, z_n)$ in which we can write
 $\p(z_1, \ldots, z_n)=z_1^2+ \cdots + z_n^2+c$.
\end{enumerate}
\end{definition}

The main result of \cite{Do00} is

\begin{theorem} [Theorem 2 in \cite{Do00}] \label{thm:exist}
 Let $(M,\omega)$ be a symplectic manifold such that the cohomology
 class $[\omega]/2\pi$ admits an integer lift $h$ in $H^2(M,\Z)$. Then
 there exists a symplectic Lefschetz pencil whose fibers
 are homologous to the Poincar\'{e} dual of $k h$, for
 $k$ large enough.
\end{theorem}

The Lefschetz pencils obtained by Theorem \ref{thm:exist} will be
called Donaldson's Lefschetz pencils. For $k$ fixed large enough,
they lie in a distinguished isotopy class of Lefschetz pencils.
Moreover they enjoy various remarkable topological and geometric
properties (see e.g.\ \cite{AMP02}). Our main result relates the
geometry of these pencils to that of Lagrangian submanifolds.
More precisely, the result that we want to prove is

\begin{theorem} \label{thm:main}
Let $(M,\omega)$ be an integral symplectic manifold. Let $\SL$ be
a Lagrangian submanifold of $M$ and let $f:\SL \to [0,1]$ be a
Morse function. Then there exists a sequence of Donaldson's
Lefschetz pencils $\phi_k:M-B_k \to \CP^1$ such that, for $k$
large, $\phi_k(\SL)$ is a smooth embedded arc in $\CP^1$.
Moreover, there exists a parametrization of this arc
$\gamma_k:[0,1] \to \phi_k(\SL)$ in such a way that
$(\gamma_k^{-1} \circ \phi_k)|_{\SL}$ is a Morse function isotopic
to $f$.
\end{theorem}

Recall that we mean that two Morse functions are isotopic if
they are isotopic {\em among Morse functions\/}, i.e.,\ there
exists a $1$-parametric family of (non-degenerate) Morse functions
connecting them.

Statements similar to Theorem \ref{thm:main} have been part of the
mathematical folklore for the past few years, after Donaldson
suggested such a picture; however, to our knowledge the details of
the statement had not been worked out, and our proof shows that
some rather unexpected technical complications do occur. It is
also worth mentioning the relationship between our result and the
work of Seidel, who showed that if two vanishing cycles of a
Lefschetz pencil can be joined by a ``matching path'' (see Section
\ref{sec:matching}), then the total space of the pencil contains a
Lagrangian sphere fibered over an arc in $\CP^1$; in the case
where $\SL\iso S^n$ and $f$ is a Morse function with only two
critical points, our result can therefore be thought of as a
converse to Seidel's construction.

It is possible to adapt our result to $S^1$-valued Morse functions.
The proof will follow exactly the same pattern and we do not
give the details but leave them to the careful reader.

The task of the next part of the article is to construct an equivalence
of sets between Hamiltonian isotopy classes of Lagrangian spheres and
the set of matching paths modulo a natural action of the fundamental
group of the space of Lefschetz pencils. The precise result will be stated in Section
\ref{sec:isotopies}. For this we will use a parametric version of Theorem
\ref{thm:main}:
\begin{theorem} \label{thm:main_par}
Let $(M,\omega)$ be an integral symplectic manifold. Let $\{ \SL_t
\}$ be a $1$-parametric family of simply connected Lagrangian
submanifolds of $M$, and let $\phi_{k,0}$ and $\phi_{k,1}$ be two
sequences of Donaldson pencils obtained using the construction of
Theorem \ref{thm:main}, adapted to the submanifolds $\SL_0$ and
$\SL_1$ and two Morse functions $f_j:\SL_j\to [0,1]$ $(j=0,1)$.
Assume moreover that $f_0$ and $f_1$ are isotopic through a family
of Morse functions $f_t$. Then there exists a sequence of families
of Donaldson pencils $\phi_{k,t}$ adapted to $\SL_t$ and $f_t$ and
coinciding with $\phi_{k,0}$ and $\phi_{k,1}$ at the ends.
\end{theorem}
The proof will be an extension of the non-parametric case; both cases rely
on a local Lemma, extending previous results of Donaldson, that will
be proved in Section \ref{sec:parametric}.


Finally, in Section \ref{sec:automorphisms} we define the group of
automorphisms $\Gamma(\phi)$ of a Lefschetz pencil and discuss its
properties. In particular, we show the asymptotic surjectivity as
$k\to\infty$ of a natural homomorphism $\rho:\Gamma(\phi)\to
\pi_0\mathrm{Symp}(M,\omega)$, and exhibit various natural
elements of its kernel. We also discuss the implications for
matching paths and the relation with the projective dual of the
discriminant curve of a projection to $\CP^2$.

\noindent {\em Acknowledgments.\/} We are grateful to Paul Seidel
for telling us a proof of Lemma \ref{lemm:match}. First author partially
supported by NSF grant DMS-0244844. Third author supported by a Post-doctoral
Fellowship from Ministerio de Educaci\'{o}n y Cultura of Spain.
Second and third authors partially supported by project
BFM2000-0024 from Ministerio de Educaci\'{o}n, Ciencia y
Tecnolog\'{\i}a of Spain. This work has been partially supported
by the European Contract Human Potential Programme, Research
Training Network HPRN-CT-2000-00101.

\section{Asymptotically holomorphic tools} \label{sec:ah-tools}

Let $(M,\omega)$ be a symplectic manifold such that $[\omega/2\pi]
\in H^2(M;\R)$ admits a lift to an integer cohomology class and
let $h\in H^2(M;\Z)$ be one such lifting. In this case there is an
hermitian line bundle $L$ with first Chern class $c_1(L)=h$, and
we can equip $L$ with a hermitian connection $\nabla$ of curvature
$-i\omega$.

Let $J$ be an almost complex structure on $M$ compatible with
$\omega$ and let $g(u,v)=\omega (u, Jv)$ be the associated metric.
We define the sequence of rescaled metrics $g_k=kg$ with
associated distance functions $d_k$. We give the following

\begin{definition} \label{def:ah}
Let $E$ be a hermitian line bundle with connection on $M$. A
sequence of sections $s_k$ of $E\otimes L^{\otimes k}$ is called
asymptotically $J$-holomorphic if it satisfies the bounds
$|\nabla^r s_k|_{g_k}=O(1)$ for $0 \leq r \leq 3$ and
$|\nabla^{r-1} \bar{\partial} s_k|_{g_k} =O(k^{-1/2})$ for $1 \leq
r \leq 3$.
\end{definition}

In this definition and throughout the text, the notation $O(1)$ means
that there exists a bound by a uniform constant depending
neither on the point of $M$ nor on the value of $k$.

If $\SL$ is a Lagrangian submanifold of $M$, then
$[\omega]|_{\SL}=0$ in $H^2(\SL,\R)$. Therefore $h|_{\SL}$ is
torsion in $H^2(\SL, \Z)$ and there exists a positive integer $p$
such that $p\, h|_{\SL}=0$. We substitute the symplectic form
$p\,\omega$ for $\omega$ on $M$, so that the first Chern class
$c_1(L)$ becomes zero when restricted to $\SL$. Therefore the line
bundle $L|_{\SL}$ is topologically trivial, and the connection
$\nabla|_{\SL}$ is flat. Since the holonomy of $\nabla|_{\SL}$
need not be trivial, we cannot expect to find a parallel section
of $L^{\otimes k}|_{\SL}$. However, by choosing a suitable
trivialization of $L^{\otimes k}|_{\SL}$, we can ensure that the
connection 1-form is bounded by a constant (w.r.t.\ the metric
$g$), which gives the following result (Lemma 2 of \cite{Au01}):

\begin{lemma} \label{lem:trivialize}
The restriction of $L^{\otimes k}$ to $\SL$ admits a section
$\s_{\SL,k}$ such that $|\s_{\SL,k}|=1$ and
$|\nabla\s_{\SL,k}|_{g_k}=O(k^{-1/2})$ at every point of $\SL$.
\end{lemma}

Moreover, if we are given a fixed contractible open subset $U\subset\SL$,
then we can additionally assume that $\nabla\s_{\SL,k}=0$ at every point of $U$.

We will also need the following lemma:

\begin{lemma} \label{lem:darboux}
Let $x\in M$. Then there exist Darboux coordinates with respect to
the symplectic form $k\omega$, $\Psi_{k,x}: B_{g_k}(x, ck^{1/2})
\to \C^n$ $($where $c>0$ is a fixed constant independent of $x$
and $k)$ such that: $(a)$ $\Psi_{k,x}(x)=0$; $(b)$
$(\Psi_{k,x})_*$ identifies the complex structure $J_x$ of $T_xM$
with the standard complex structure $J_0$ on $\C^n$; and $(c)$
$\Psi_{k,x}$ is approximately isometric, i.e.,\ the map
$\Psi_{k,x}$ satisfies $|\nabla^r \Psi_{k,x}|=O(1)$ for $r=1,2,3$
and $|\nabla^{r-1}\bbd \Psi_{k,x}(z)|= O(k^{-1/2}|z|)$ for
$r=1,2,3$. Moreover, if $x\in \SL$ we can take $\Psi_{k,x}$ to map
the Lagrangian submanifold to $\R^n$.
\end{lemma}

\begin{proof}
 The first part of this result is Lemma 3 in \cite{Au00} (see also
\cite{Do96}). In the
 case where $x$ lies in the Lagrangian submanifold $\SL$ of $M$, we
 use Weinstein neighborhood theorem to make the charts
 $\Psi_{k,x}$ map $\SL$ to $\R^n$.
\end{proof}

Note that the estimates in Lemmas \ref{lem:trivialize} and \ref{lem:darboux}
depend on the geometry of the submanifold $\SL$ (in particular its
injectivity radius with respect to the metric $g$).

\begin{definition} \label{def:real}
An asymptotically holomorphic sequence $s_k$ of sections of
$L^{\otimes k}$ on $M$ will be called real if
 $$
 s_k|_{\SL}= f_k \cdot \s_{\SL,k},
 $$
for some non-negative real function $f_k$.
\end{definition}

For example, the sections constructed in Lemma \ref{lemma:local} below are
real. We are interested in this property because we will show later how to
preserve it under the usual local perturbations in asymptotically
holomorphic theory. In particular, we use from \cite{Au01} the following

\begin{theorem} \label{thm:auroux}
Given a symplectic manifold $(M,\omega)$ and a Lagrangian
submanifold $\SL$, there exists an asymptotically holomorphic
sequence of sections $s_k:M \to L^{\otimes k}$ such that
\begin{enumerate}
 \item $\forall x\in \SL$, $|s_k(x)|\geq \eta$,
 for some $\eta>0$ independent of $k$.
 \item $\forall x\in M$, $|s_k(x)|\leq e^{-d_{k}(x,\SL)^2/5}$.
 \item $\forall x\in M$, such that $d_{k}(x,\SL)\geq 2k^{1/6}$,
 we have $s_k(x)=0$.
 \item $s_k$ is real.
\end{enumerate}
\end{theorem}

In order to prove this Theorem we need the following
\begin{lemma} \label{lemma:local}
Let $x\in M$. There exists a sequence of asymptotically holomorphic sections
$\sref{x}$ satisfying the following estimates
\begin{itemize}
 \item $|\sref{x}(y)| \geq \frac{1}{2}$ for all $y\in B_{g_k}(x,1)$.
 \item $|\sref{x}|_{C^2} \leq p(d_k(x,y)) e^{-d_k(x,y)^2/5}$,
 where $p(t)$ is a fixed real polynomial.
 \item $\sref{x}$ has support in the ball $B_{g_k}(x, 2k^{1/6})$.
 \item If $x\in \SL$, there is a real non-negative function $f_k:\SL\to \R$
 such that $\sref{x}|_{\SL} = f_k \cdot \s_{\SL,k}$.
\end{itemize}
\end{lemma}

\begin{proof}
The existence of reference sections $\sref{x}$ satisfying the
first three properties follows from Lemma 3 in \cite{Au97} (see
also Proposition 11 in \cite{Do96}). To get the last property,
consider the Darboux chart $\Psi_{k,x}: B_{g_k}(x,ck^{1/2})\to
\C^n$ provided by Lemma \ref{lem:darboux}, which maps the
Lagrangian submanifold to $\R^n$. Now we trivialize the positive
line bundle $L^{\otimes k}$ in $B_{g_k}(x,ck^{1/2})$ following
radial directions. This trivialization yields a radially parallel
local section $s_r^k$, such that
$|\nabla(\sigma_{\SL,k}/(s_r^k|_{\SL}))|= O(k^{-1/2})$ by Lemma
\ref{lem:trivialize}. The localized sections obtained in Lemma 3
in \cite{Au97} are of the form
 \begin{equation} \label{eqn:sref}
 \sref{x,0}= e^{-|z|^2/4} \chi(|z|) s_r^k,
 \end{equation}
where $\chi(x)$ is a real positive cut-off function which equals $1$ in a
ball of $g_k$-radius $k^{1/6}$ and $0$ outside of a ball of
$g_k$-radius $2k^{1/6}$. The sections $\sref{x,0}$ satisfy all desired
properties except that over $\SL$ they are real multiples of $s_r^k$ instead
of $\s_{\SL,k}$. However, we can choose a function $\phi:B_{g_k}(x,2k^{1/6})
\to\R$ such that $(i)$ $|\phi|_{C^2}=O(k^{-1/2})$ and $(ii)$ $\s_{\SL,k}=e^{i\phi}
s_r^k|_{\SL}$. Then $(i)$ implies that $\sref{x}=e^{i\phi}\sref{x,0}$ still
satisfies the first three properties, while $(ii)$ implies that $\sref{x}$
is a positive multiple of $\sigma_{\SL,k}$ over $\SL$.
\end{proof}

\begin{remark}\label{rmk:holomorphic}
Assume that there exists a neighborhood $U$ of $x$, containing
$B_{g_k}(x,2k^{1/6})$, over which
$J$ is integrable and standard in local Darboux coordinates. If $x\in \SL$,
assume moreover that $\s_{\SL,k}$ is covariantly constant over $U\cap \SL$.
Then $\sref{x}$ is $J$-holomorphic over $B_{g_k}(x,k^{1/6})$.
\end{remark}

Let us now sketch a proof of Theorem \ref{thm:auroux}.

\begin{proof}[Proof of Theorem \ref{thm:auroux}.] Take a set of
points $S_k$ of $\SL$ such that the balls $B_{g_k}(p,1)$, $p\in
S_k$, cover $\SL$ and any point of $\SL$ belongs to at most $N$
such balls (where $N$ is a constant depending only on the
dimension). Then we define a global section
 $$
 s_k= \sum_{x_j\in S_k} w_j \sref{x_j},
 $$
where $\sref{x_j}$ are provided by Lemma \ref{lemma:local} and $w_j\in\C$ are
constants. If we
choose $w_j= 1$, then the sequence is asymptotically holomorphic
and satisfies all the properties listed in the statement of the
Theorem. \end{proof}

We also have the usual notion of estimated transversality
\cite{Au97}:

\begin{definition} \label{def:eta-trans}
 A section $s_k$ of the bundle $E\otimes L^{\otimes k}$ is $\eta$-transverse
 to $0$ if, for every $x\in M$ such that
 $|s_k(x)|<\eta$, $\nabla s_k(x)$ has
 a right inverse $\theta_k$ such that $|\theta_k|_{g_k}<\eta^{-1}$.
\end{definition}

\begin{lemma} \label{lem:lag-trans}
 Let $\SL$ be a Lagrangian submanifold of $M$ and let $\varphi_k$ be an
 asymptotically holomorphic sequence of complex valued functions
 defined in some neighborhoods of $\SL$. Suppose that $h_k= \varphi_k|_{\SL}$
 takes values in an immersed curve $C\subset \C$. If $h_k$
 satisfies the following property (identifying $\nabla\nabla h_k(x)$ with
 a linear map $T_x\SL\to T_x\SL \iso T_x \SL\otimes T_{h_k(x)}C \subset T_x\otimes\C$
 by means of the metric $g_k$):
 $$
  |\nabla h_k(x)| <\eta \Rightarrow \nabla\nabla h_k(x): T_x\SL\to
  T_x\SL\text{ multiplies the length of vectors by at least }\eta
 $$
 for any $x\in \SL$, then $\nabla \varphi_k$ is $\eta/2$-transverse to
 $0$ on a $g_k$-neighborhood of $\SL$ of uniform radius.
\end{lemma}

\begin{proof}
 For $x\in \SL$, if $|\nabla \varphi_k(x)|<\eta$ then $|\nabla
 h_k(x)|<\eta$, hence $\nabla\nabla h_k(x): T_x\SL\to
 T_x\SL$ multiplies the length of vectors
 by at least $\eta$. But $T_xM = T_x\SL \oplus J(T_x\SL)$ and the
 direct sum is an orthogonal one. Since the Hessian $\nabla\nabla
 \varphi_k(x)$ is approximately holomorphic, it equals the
 complexification of $\nabla\nabla h_k(x)$ up to some error of
 order $O(k^{-1/2})$. Therefore, $\nabla \nabla \varphi_k(x)$ multiplies the
 length of vectors by at least $3\eta/4$. It can be checked that this
 property is equivalent to the $3\eta/4$-transversality of $\nabla\varphi_k$
 in the sense of Definition \ref{def:eta-trans}. From this it
 follows the $\eta/2$-transversality in a $g_k$-neighborhood of
 $\SL$ of uniform radius.
\end{proof}

\begin{definition}\label{def:lag-trans}
 Let $h_k:\SL \to \C$ be a sequence of functions taking values
 in an immersed curve of $\C$. We say
 that $\nabla h_k$ is $\eta$-transverse if, for every $x\in \SL$ such that
 $|\nabla h_k(x)|<\eta$, the map $\nabla\nabla h_k(x): T_x\SL\to
 T_x\SL$ multiplies the length of vectors by at least $\eta$.
\end{definition}

Finally, recall the following result, due to Donaldson
\cite{Do00}, which implies Theorem \ref{thm:exist}:

\begin{proposition} \label{prop:get_tr}
Given an asymptotically holomorphic sequence of sections
$s_{k}^1\oplus s_{k}^2$ of $L^{\otimes k}\oplus L^{\otimes k}$,
and given $\d>0$, there exists a sequence of sections
$\s_{k}^1\oplus \s_{k}^2$ with $|\s_{k}^j-s_{k}^j|_{C^2}\leq \d$,
for $j=1,2$, satisfying the following properties:
 \begin{enumerate}
 \item $\s_{k}^1$ is $\e$-transverse to zero over $M$, for some
 uniform $\e>0$.
 \item $\s_{k}^1\oplus\s_{k}^2$ is $\e$-transverse to zero over $M$.
 \item Denoting by $Z_{k,\e}=\{ p\in M: |\s_{k}^1|\leq \e \}$, the
 map $\bd \left( \s_{k}^2 / \s_{k}^1\right)$ is $\e$-transverse to
 zero in $M-Z_{k,\e}$.
 \end{enumerate}
 Moreover, after a small perturbation around the critical points we
 can ensure that $\p_k= \s_k^2/\s_k^1$ is a symplectic Lefschetz pencil.
\end{proposition}

\section{Deformations of a Morse function} \label{sec:morse}

The first step in our proof of Theorem \ref{thm:main} is to
perturb the given Morse function $f$ on $\SL$ to a suitable
sequence of functions $h_k$. We need this sequence to satisfy that
$\nabla h_k$ be transverse to zero in an estimated way. This is a
necessary condition to have the functions $h_k$ arise as restrictions of Donaldson
pencils (by a converse to Lemma \ref{lem:lag-trans}).
There are many ways of achieving such transversality,
however we have to be careful because we will try to approximate $h_k$
by a combination of sections constructed by Lemma
\ref{lemma:local} that needs to be controlled in a precise way. This is
the reason of the extra requirements that we will impose.

Let $f$ be a Morse function on $\SL$. Consider a critical point
$p$ of $f$, i.e.,\ $\nabla f(p)=0$. Then for a suitable chart
$\phi=(x_1,\ldots, x_n): U\to \R^n$ with $\phi(p)=0$ and which is
isometric at $p$, we can write $f=c+ \l_1 x_1^2 +\cdots +\l_n
x_n^2$ for suitable non-zero real numbers $\lambda_1, \ldots,
\lambda_n$. We start by deforming $f$ in a neighborhood of each
critical point so that we can write $f=c+ \e_1 x_1^2+\cdots + \e_n
x_n^2$, where $\e_i=\pm 1$ depending on whether $\l_i$ is positive
or negative. We call $f$ again this new function. Moreover, for
simplicity we also deform the almost complex structure $J$ in a
neighborhood of every critical point so that it is integrable and
standard in Darboux coordinates. Therefore, in such a
neighborhood, $J$ and $g$ have a standard form. Let $r_1>0$ be a
constant such that $f$ and $J$ can be written in such a standard
way over the entire ball of radius $r_1$ around each critical
point.

By Lemma \ref{lem:darboux} there exists some $r_2>0$ such that all
Darboux charts $\Psi_{k,x}$ are ``approximately isometric'' in
$B_{g}(x,r_2)$, i.e.,\ the differential of $\Psi_{k,x}$ and its
inverse distort the metric in a controlled manner at any point of the
chart (in the sense of Lemma \ref{lem:darboux}). Let
$c_0$ be the minimum of $r_1$ and $r_2$. Note that for a critical
point $p$, $\Psi_{k,p}$ is actually holomorphic and an isometry,
since $J$ and $g$ are standard in $B_{g_k}(p, c_0k^{1/2})$.

We need the following technical Lemma that will be used to define
deformations of a Morse function.

\begin{lemma} \label{lemma:scale}
Let $l:\R_+ \to \R_+$ be a smooth function with $l'(t)<0$ and
$\displaystyle{\frac{l'(t)}{l(t)} \geq -\frac{3}{4t}}$ for all
$t>0$. Then the smooth function $ l(|x|) \cdot x: \R^n-\{0\} \to
\R^n- \{0\}$ satisfies the following properties:
\begin{enumerate}
\item $d ( l(|x|)\cdot x ) = l(|x|)\Big( (dx_1, \ldots, dx_n) +
\displaystyle{\frac{l'(|x|)}{l(|x|)\, |x|}} (x_1x_1 dx_1+ \cdots +
x_1x_n dx_n, \ldots, x_nx_1dx_1 +\cdots + x_nx_ndx_n) \Big)$,
i.e., using matrix notation, $\displaystyle d(l(|x|)\cdot
x)=l(|x|)\,\mathrm{Id}+\frac{l'(|x|)}{|x|}\,(x\otimes x^t)$. \item
$\det(d ( l(|x|)\cdot x ) ) = l(|x|)^n\left(
\displaystyle{1+\frac{|x| l'(|x|)}{l(|x|)}}\right) \geq
l(|x|)^n/4$. \item $|d (l(|x|)\cdot x )| \leq l(|x|)$. \item The
minimum eigenvalue of $d ( l(|x|)\cdot x)$ is greater or equal to
$l(|x|)/4$. \item $l(|x|)\cdot x$ is a global diffeomorphism.
\end{enumerate}
\end{lemma}

\begin{proof}
The first property is a simple computation. The next three properties
follow from the expression for the linear map $d(l(|x|)\cdot x)$, which can
be diagonalized by means of an orthogonal transformation mapping $x$ to the
first coordinate axis. It is then clear that the eigenvalues of
$d(l(|x|)\cdot x)$ are $l(|x|)$ with multiplicity $n-1$, and
$l(|x|)+l'(|x|)|x|$ with multiplicity $1$. The last property follows
from the previous considerations.
\end{proof}

\begin{proposition} \label{prop:deform}
Let $f:\SL \to [0,1]$ be a Morse function and let $D>0$ be a large
enough constant. Then there exists a sequence $h_k:\SL\to [0,k^{1/2}]$ of
deformations of $f$ satisfying (with respect to the metric $g_k$ over $\SL$)
 \begin{enumerate}
  \item $|\nabla h_k|=O(D)$, $|\nabla\nabla h_k|=O(1)$.
  \item $\nabla h_k$ is $\eta$-transverse to zero, for some
  $\eta>0$ independent of $k$ and $D$.
  \item $|\nabla^3 h_k|=O(1/D)$.
\end{enumerate}
\end{proposition}

To prove the main result, we will need to choose a sufficiently large
value of the constant $D$; this requires a slightly more complicated
argument (otherwise one could set $\epsilon=0$ in the formulas below).

\begin{proof}
Denote by $\{p_j\}$ the set of critical points of $f$. We are
assuming that $f$ is of the form
 \begin{equation*}
  f(x_1,\ldots, x_n)= c_j + \sum \pm x_i^2,
 \end{equation*}
in a neighborhood $B_g(p_j,c_0)$ of $p_j$. Scaling the coordinates
we get that, over $B_{g_k}(p_j, c_0k^{1/2})$, $f$ is expressed as
 $$
 f(x_1,\ldots, x_n)= c_j + \frac{1}{k}\sum \pm x_i^2 ,
 $$
in the chart $\Psi_{k,p_j}$. These new coordinates are actually
isometric over the ball $B_{g_k}(p_j,c_0k^{1/2})$, because of the
choice of almost-complex structure $J$ made at the beginning of this
section. Now define
the sequence $f_k=k^{1/2} f$. The derivatives of $f_k$ are
\begin{eqnarray}
  f_k(x_1, \ldots, x_n) & = & k^{1/2}c_j + \frac{1}{k^{1/2}}\sum
  (\pm x_i^2), \nonumber\\
  \nabla f_k(x_1, \ldots, x_n) & = & \frac{1}{k^{1/2}} \sum (\pm
  2x_idx_i), \label{eq:Morse} \\
  \nabla\nabla f_k(x_1, \ldots, x_n) & = & \frac{2}{k^{1/2}} \sum (\pm
  dx_i\otimes dx_i). \nonumber
\end{eqnarray}
So $f_k$ satisfies, without any perturbation,
 \begin{eqnarray}
 |\nabla f_k|_{g_k} & = & O(1) \nonumber \\
 |\nabla \nabla f_k|_{g_k} & =& O(k^{-1/2}). \label{eqn:cont}
 \end{eqnarray}

Now take a cut-off function $l_k:[0,k^{1/2}c_0] \to [0,k^{1/4}]$
satisfying the following properties:
 \begin{itemize}
  \item $l_k(t)=k^{1/4}$, for $t\in [0, D]$.
  \item $l_k(t)= \displaystyle{\frac{a\,k^{1/4}}{t^{1/2+\e}}}$, for
    $t\in [2D, k^{1/2}c_0/2]$, where $a$ and $\e$ are positive constants
    to be adjusted below ($a$ converges to a fixed positive value as
    $k\to\infty$, while $\e=\e(k)\to 0$ as $k\to \infty$).
  \item $l_k(t)=1$, for $t>3k^{1/2}c_0/4$.
  \item In the interval $[D, 2D]$, $l_k$ is decreasing,
    $|l_k'(t)| \leq
    \epsilon_2\,k^{1/4}/D$ and $|l_k''(t)| \leq C k^{1/4}/D^2$,
    $|l_k'''(t)| \leq C k^{1/4}/D^3$
    for constants $\epsilon_2,C$ independent of $k$.
  \item In the interval $[k^{1/2}c_0/2, 3k^{1/2}c_0/4]$, $l_k$ is
   decreasing, $|l_k'(t)| \leq \epsilon_2'\,
   k^{-1/2}$ and $|l_k''(t)| \leq C' k^{-1}$,
   $|l_k'''(t)| \leq C' k^{-3/2}$ for constants
   $\epsilon_2',C'$ independent of $k$.
  \item The following inequality holds over $(D, 3k^{1/2}c_0/4)$:
   \begin{equation}
    0>\displaystyle{\frac{l_k'(t)}{l_k(t)} \geq -\frac{1/2+\e}{t}}
    \label{eq:above}
   \end{equation}
\end{itemize}

The following figure gives the shape of $l_k(t)$:

\includegraphics{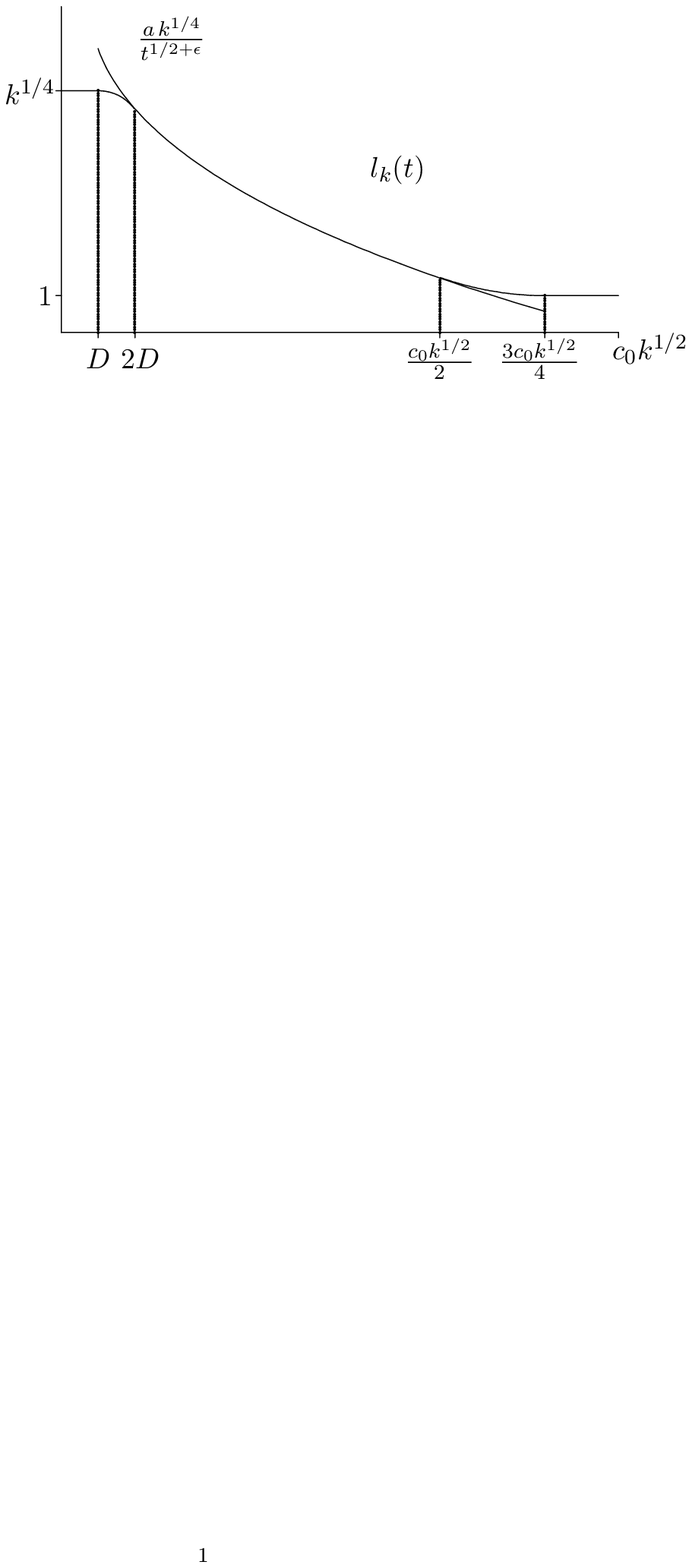}

We adjust $a$ and $\e$ in the following way: the function
${\displaystyle \frac{a\,k^{1/4}}{t^{1/2+\e}}}$ equates $k^{1/4}$
for $t=\frac32 D$ and equates $1$ for $t=0.7 c_0 k^{1/2}$.
Therefore
 $$
 \left\{ \begin{array}{ccl}
  a &=& \left({\displaystyle \frac32} D\right)^{1/2+\e} \\
  \e &=& {\displaystyle \frac{\ln  {\displaystyle\frac{3D}{1.4c_0}}}{\ln k -2 \ln
   {\displaystyle\frac{3D}{1.4c_0}}}}
 \end{array} \right.
 $$
The inequality $k\geq \left({\displaystyle\frac{3D}{1.4c_0}}\right)^6$
must hold in order to have $0<\e\leq \frac14$. With these choices of $a$
and $\epsilon$ it is easy to find a function $l_k(t)$ with the desired
properties; observe in particular that
(\ref{eq:above}) automatically holds when $l_k$ is of the form
$l_k(t)=F(t) \cdot {\displaystyle
\frac{a\, k^{1/4}}{t^{1/2+\e}}}$ for some increasing function
$F$ (by differentiating the identity $\log l_k =\log F -(\frac12+\e) \log t + C$).

Now we define the function $\hat{l}_k(x_1, \ldots, x_n)=
l_k(|(x_1, \ldots, x_n)|)$ in a small neighborhood of the critical
points. This implies the following formulae for $\hat{l}_k$:
  \begin{eqnarray} \label{eqn:hatl}
    |\nabla \hat{l}_k(x)| &\leq& C_1|l'_k| \leq C_1 \frac{|l_k|}{|x|}  \nonumber\\
    |\nabla^2 \hat{l}_k(x)| &\leq& C|l''_k|+ C\frac{|l'_k|}{|x|}
     \leq C_2 \frac{|l_k|}{|x|^2}  \\
    |\nabla^3 \hat{l}_k(x)| &\leq& C|l'''_k|+ C\frac{|l''_k|}{|x|}
     + C\frac{|l'_k|}{|x|^2}\leq C_3 \frac{|l_k|}{|x|^3} \nonumber
  \end{eqnarray}
for some uniform constants $C_1,C_2,C_3$. These are simple
computations in the middle interval and in $[D,2D]$ and
$[k^{1/2}c_0/2,3k^{1/2}c_0/4]$, using the bounds on the derivatives of
$l_k$.

We construct a perturbation of $f$ by setting:
 $$
 h_k(x)= \left\{ \begin{array}{ll} f_k( \hat{l}_k(x)\cdot x),\qquad
 & \forall x\in \bigcup\limits_{j} B_{g_k}(p_j, k^{1/2}c_0) \\
 f_k(x), & \text{otherwise.} \end{array} \right.
 $$
Outside the union of the balls $B_{g_k}(p_j, 3k^{1/2}c_0/4)$, the
estimates (\ref{eqn:cont}) imply that $h_k$ satisfies (i).
Moreover, in that region, there is a constant $\xi>0$ independent
of $k$ such that $|\nabla h_k|_{g_k} \geq \xi$, which implies that $\nabla
h_k$ is $\xi$-transverse to zero, so (ii) also holds.

In the ball $B_{g_k}(p_j, D)$, we have $h_k=k^{1/2}c_j + \sum (\pm
x_i^2)$, so $h_k$ satisfies (i) and (ii).

Now let us see what happens in the annuli $B_{g_k}(p_j,
3k^{1/2}c_0/4)- B_{g_k}(p_j, D)$. We start by computing the
derivatives of $h_k$,
  \begin{eqnarray}
  \nabla h_k(x) &=& \nabla f_k(\hat{l}_k(x)\cdot x) \circ
  (d \hat{l}_k(x) \cdot x + \hat{l}_k(x) \cdot \Id),
  \label{eqn:der_1} \\
  \nabla\nabla  h_k(x) &=& \nabla\nabla f_k(\hat{l}_k(x)\cdot x)
  \circ ((d \hat{l}_k(x)
  \cdot x + \hat{l}_k(x) \cdot \Id)\otimes (d \hat{l}_k(x) \cdot x +
  \hat{l}_k(x) \cdot \Id)) + \nonumber\\ & & + \nabla f_k(\hat{l}_k(x)\cdot x)
  (\nabla\nabla \hat{l}_k(x) \cdot x +  d\hat{l}_k(x) \otimes \Id+\Id
  \otimes d\hat{l}_k(x) ). \label{eqn:der_2}
  \end{eqnarray}
Note that we have, using (\ref{eq:Morse}),
 $$
 |\nabla f_k(\hat{l}_k(x)\cdot x)|= \frac{2}{k^{1/2}}
 |\hat{l}_k(x) \cdot x|=\frac{2}{k^{1/2}} l_k(|x|)|x|.
 $$
Also, consider $d (\hat{l}_k(x)\cdot x)=d \hat{l}_k(x) \cdot x +
\hat{l}_k(x) \cdot \Id$. The function $l_k$ satisfies that
$l_k'(t)<0$ and $\displaystyle{\frac{l_k'(t)}{l_k(t)} \geq
-\frac{3}{4t}}$. Hence we may apply Lemma \ref{lemma:scale} to get
that the eigenvalues of $d (\hat{l}_k(x)\cdot x)$ are in the range
$[\frac{1}{4}\,l_k(|x|), l_k(|x|)]$.

By construction, we have, in the region $D \leq t \leq
3c_0k^{1/2}/4$, that
 $$
 \frac{c_1 D^{1/2+\e} k^{1/4}}{t^{1/2+\e}} \leq l_k(t)\leq \frac{c_2
 D^{1/2+\e}k^{1/4}}{t^{1/2+\e}},
 $$
for some positive constants $c_1,c_2$ independent of $k$ and $D$.
Therefore (\ref{eqn:der_1}) implies
  $$
  |\nabla h_k(x)| \leq \frac{2}{k^{1/2}} l_k(|x|) |x|
  \, l_k(|x|) \leq  \frac{C D^{1+2\e}}{|x|^{2\e}} \leq O(D),
  $$
and analogously,
  $$
  |\nabla h_k(x)|\geq \frac{c D^{1+2\e}}{|x|^{2\e}} \geq
  \frac{c D^{1+2\e}}{(c_0k^{1/2})^{2\e}} \geq \xi
  $$
for some uniform $\xi$, thus proving (ii). To finish proving (i), we have
to bound the second derivative.
  \begin{eqnarray*}
    |\nabla\nabla  h_k(x)| &=& \big| \nabla\nabla f_k(\hat{l}_k(x)\cdot x)
  \circ ((d \hat{l}_k(x)
  \cdot x + \hat{l}_k(x) \cdot \Id)\otimes (d \hat{l}_k(x) \cdot x +
  \hat{l}_k(x) \cdot \Id)) + \\ & & + \nabla f_k(\hat{l}_k(x)\cdot x)
  (\nabla\nabla \hat{l}_k(x) \cdot x +  d\hat{l}_k(x) \otimes \Id+\Id
  \otimes d\hat{l}_k(x) ) \big| \\
  &\leq & \frac{2}{k^{1/2}} l_k( |x|)^2 +
  \frac{2}{k^{1/2}} l_k(|x|) |x| \Bigl[C_2\,\frac{l_k(|x|)}{|x|} +2|l_k'(|x|)|
  \Bigr] \\ &\leq &C
  \frac{l_k(|x|)^2}{k^{1/2}} \leq O(1),
\end{eqnarray*}
using (\ref{eqn:hatl}) and $|l_k'| |x| \leq \frac34 \, l_k$.

Finally we prove (iii) by considering the third derivative of $h_k$.
If the considered point is not close to a critical point then $f_k$ has
not been perturbed, i.e.,\ $h_k=f_k$ and we have
$|\nabla\nabla\nabla f_k|_{g_k}= O(k^{-1})$. To bound
$\nabla\nabla\nabla h_k$ in a
ball $B_{g_k}(p_j, c_0 k^{1/2})$ around a critical point $p_j$, we
compute from (\ref{eqn:der_2}),
 \begin{eqnarray*}
   |\nabla\nabla\nabla  h_k |& = & \Big| \sum\nolimits'
   \nabla\nabla f_k(\hat{l}_k(x)\cdot x) \circ
   \Big( (d\hat{l}_k(x)\cdot x +\hat{l}_k(x)\cdot \Id) \otimes
   (\nabla\nabla\hat{l}_k(x)\cdot x +  d\hat{l}_k(x) \otimes\Id +\\
   & & + \Id \otimes d\hat{l}_k(x))\Big)
   + \nabla f_k(\hat{l}_k(x) \cdot x) \circ
   \left(\nabla\nabla\nabla  \hat{l}_k(x) \cdot x + \sum\nolimits'
   \nabla\nabla\hat{l}_k(x) \otimes \Id\right) \Big| \\
   & \leq & 3\frac{2}{k^{1/2}}\,l_k(|x|) (C_2+2) \frac{l_k(|x|)}{|x|} +
   \frac{2 |x|l_k(|x|)}{k^{1/2}} (C_3+3C_2)\frac{l_k(|x|)}{|x|^2}\\
   & \leq & C \frac{l_k(|x|)^2}{k^{1/2}|x|} = O \left({\frac1D}\right).
 \end{eqnarray*}
where $\sum'$ denotes the sum of three terms obtained by cyclic
permutation of the three factors in $T^*M\otimes T^*M\otimes
T^*M$. This estimate completes the proof of the Proposition.
\end{proof}

The purpose of the bound on $|\nabla^3 h_k|$ given in Proposition
\ref{prop:deform} is to let us control the error produced by replacing
$h_k$ by its order two Taylor approximation at a given point.
Fix any point $x\in \SL$ and consider the approximately isometric
chart $\Psi_{k,x}$ of Lemma \ref{lem:darboux}. Let $P_{2,x}$ be
the degree $2$ Taylor polynomial of $h_k$ at $x$ in these coordinates.
Then we have in particular
$|h_k(y)- P_{2,x}(y)| \leq  |\nabla\nabla\nabla h_k(\xi)| |y|^3$,
where $\xi$ is a point of the segment joining $0$ and $y$; hence, it
follows from Proposition \ref{prop:deform} (iii) that
$$|\nabla^i (h_k(y)- P_{2,x}(y))| \leq {\frac{C}{D}} |y-x|^{3-i},
\ \mathrm{for}\ y\in \SL\cap B_{g_k}(x,k^{1/2}c_0),\ i=0,1,2.$$

\section{Transversality in a neighborhood of the Lagrangian.}
\label{sec:trans_1}

The objective now is to define two asymptotically holomorphic
sections $s_k^1$ and $s_k^2$ whose quotient approximates $h_k$.
This cannot be done directly because $|h_k|$ is not $O(1)$ and
therefore if we take the sequence $s_k^1$ as given in Theorem
\ref{thm:auroux}, then $s_k^2=h_k \cdot s_k^1$ will not be
bounded. We shall work in a different way.

\begin{lemma} \label{lemma:pair}
Let $h_k$ be a sequence of functions satisfying the bounds of
Proposition \ref{prop:deform}. Define a new sequence
  \begin{eqnarray*}
   (h_k^1, h_k^2): \SL & \to & S^1 \\
    x & \to & (\cos(h_k(x)), \sin(h_k(x))).
  \end{eqnarray*}
Then these sections satisfy:
  \begin{enumerate}
    \item $|h_k^j|=O(1)$, and $|\nabla^r h_k^j|=O(D)$, $r=1,2$, $j=1,2$.
    \item $\nabla (h_k^1, h_k^2)$ is $\eta$-transverse to zero,
    for some $\eta>0$ independent of $k$.
    \item Fix any point $x\in \SL$ and consider the chart
  $\Psi_{k,x}$ of Lemma \ref{lem:darboux}. Let $P_{2,x}$ be
  the degree two Taylor polynomial approximating $h_k$ at $x$. Then
   \begin{itemize}
  \item $|\nabla^i (h_k^1(y)- \cos (P_{2,x}(y)))|\leq
  \frac{C}{D}|y-x|^{3-i}$,
  for $y\in \SL\cap B_{g_k}(x,k^{1/2}c_0)$, $i=0,1,2$.
  \item $|\nabla^i (h_k^2(y)- \sin (P_{2,x}(y)))|\leq
  \frac{C}{D}|y-x|^{3-i}$,
  for $y\in \SL\cap B_{g_k}(x,k^{1/2}c_0)$, $i=0,1,2$.
    \end{itemize}
  \end{enumerate}
\end{lemma}

\begin{proof}
The three properties are simple to check by using the chain rule and the
Taylor approximation theorem.
\end{proof}

\begin{proposition} \label{prop:nice_funct}
Fix a sequence of functions $h_k:\SL \to \R$ satisfying the bounds
of Proposition \ref{prop:deform}. Then for any $\delta>0$ there
are two real asymptotically holomorphic sections $s_k^1$ and
$s_k^2$ of $L^{\otimes k}$, supported over a fixed neighborhood of
the zero section in $T^*\SL$ and such that
 $$
 \left|(\cos(2h_k), \sin (2h_k)) - \left(
 \frac{s_k^1+is_k^2}{s_k^1-is_k^2} \right) \big|_{\SL}\right|_{C^2}
 \leq \delta.
 $$
\end{proposition}

The reason why we need a $C^2$ approximation is to make sure that the
restriction to $\SL$ of the pencil that we construct remains a Morse
function isotopic to $h_k$.

\begin{proof}
Construct the associated sequence $(h_k^1, h_k^2)=(\cos
(h_k),\sin(h_k))$ which satisfies the bounds of Lemma
\ref{lemma:pair}. Assume for a moment that we have two
asymptotically holomorphic sections $s_k^0$ and $s_k^1$ defined
over $M$ whose quotient approximates $h_k^1$ with error
 \begin{equation}
 \left|h_k^1- \left( \frac{s_k^1}{s_k^0} \right)\big|_\SL\right|_{C^2}
 \leq \gamma. \label{eq:err1}
 \end{equation}
Analogously, assume that we have $s_k^0$ and $s_k^2$ whose
quotient approximates $h_k^2$ with error
 \begin{equation}
 \left|h_k^2- \left( \frac{s_k^2}{s_k^0} \right)\big|_\SL\right|_{C^2}
 \leq \gamma. \label{eq:err2}
 \end{equation}
Assume moreover that $( s_k^2/s_k^1 )\big|_{\SL} \in \RP^1=\R\cup\{\infty\}$.
Now note that
 $$
 \frac{h_k^1+i h_k^2}{h_k^1- i h_k^2}= (h_k^1+ i
 h_k^2)^2=(\cos(2h_k), \sin (2h_k))
 $$
and therefore
 $$
 \left|(\cos(2h_k), \sin (2h_k)) - \left(  \frac{s_k^1+i s_k^2}{s_k^1- i
 s_k^2} \right)\big|_{\SL}\right|_{C^2} \leq \delta,
 $$
where $\gamma$ depends on $\delta$. Note that since
$(s_k^2/s_k^1)|_{\SL} \in \RP^1$, the term $((s_k^1+i
s_k^2)/(s_k^1- i s_k^2))|_{\SL}$ takes values on the unit circle.

To conclude we just have to find asymptotically holomorphic real
sequences $s_k^0$, $s_k^1$ and $s_k^2$ satisfying equations
(\ref{eq:err1}) and (\ref{eq:err2}). The condition of being real
will ensure that $(s_k^2 /s_k^1)|_{\SL} \in \RP^1$. Let us treat
the case of $h_k^1$, the other one being analogous.

Let $s_k^0$ be a sequence of asymptotically holomorphic real
sections as given by Theorem \ref{thm:auroux}. Recall from Lemma
\ref{lem:trivialize} that we have a trivializing section
$\sigma_{\SL,k}$ of $L^{\otimes k}|_{\SL}$. There exists a set of
points $\{x_j\}_{j\in J}$ in $\SL$ such that the sequence of
sections is
 $$
 s_k^0= \sum_{j\in J} w_j \sref{x_j},
 $$
where $\sref{x_j}$ are the local sections given by Lemma
\ref{lemma:local} and $w_j=1$. Remember that the points $x_j$
are chosen in such a way
that the balls $B_{g_k}(x_j,1)$ cover $\SL$ and any point of $\SL$
belongs to at most a fixed number of such balls.
Recall that $s_k^0$ is real, and that
$\left|s_k^0|_{\SL}\right| \geq \eta$, for some constant $\eta>0$
independent of $k$.

For any given point $x_j$, take an approximately holomorphic chart
$\Psi_{k,x_j}=(z_1, \ldots, z_n)$ given by Lemma
\ref{lem:darboux}. This chart maps the Lagrangian submanifold to
$\R^n \subset \R^n \oplus i \R^n = \C^n$. We consider the second
order Taylor expansion of $h_k$ at $p_j$, which we denote by
$P_{2,j}$. This Taylor expansion can be understood as a real
degree $2$ polynomial in the real parts of the Darboux
coordinates, but $P_{2,j}$ can also be viewed as a complex
polynomial (with real coefficients) in the whole $\C^n$. Recall
that by Lemma \ref{lemma:pair} we have
 \begin{equation}
 |\nabla^i (h_k^1(y)- \cos(P_{2,j}(y)))|
  \leq \frac{C}{D}|y-x_j|^{3-i}.
 \label{eq:good_approx}
 \end{equation}

We define the section
 $$
  s_k^1 = \sum_{x_j \in S} w_j \cos(P_{2,j}) \, \sref{x_j},
 $$
and use $s_k^1/s_k^0$ as a candidate to approximate our sequence
of functions $h_k^1$. It is easy to check that $s_k^1$ is real
(since the sections $\sref{x_j}$ and the polynomials $P_{2,j}$ are
real); however the existence of a uniform bound on $s_k^1$ depends
on the growth of the function $\cos(P_{2,j})$ away from $x_j$. For
this purpose, note that we may assume from the beginning that
$|\nabla\nabla h_k| < \frac{1}{10}$, and therefore that
$\cos(P_{2,j})$ has a quadratic exponential term of order less
than $\frac{1}{8}$. Indeed, the estimate of Proposition
\ref{prop:deform} ensures that such a bound can easily be achieved
by multiplying $h_k$ by a small positive real constant
(independent of $k$ and $D$). In this case, it is very simple to
check that $\cos (P_{2,j}) \sref{x_j}$ is asymptotically
holomorphic and still presents an exponential decay (just with a
lower constant). Notice that this loss in the exponential decay
only happens close to the critical points, since away from them
the degree $2$ term of $P_{2,j}$ is $O(k^{-1/2})$ and so for large
$k$ the extra factor $\cos(P_{2,j})$ does not affect the decay
properties of the product.

The $C^0$ error made with the approximation $s_k^1/s_k^0$ is
 \begin{eqnarray}\label{eq:error_0}
 \left|h_k^1- \left(\frac{s_k^1}{s_k^0}\right)\big|_{\SL}\right|
  &=& \frac{1}{|s_k^0|} |h_k^1
 \cdot s_k^0 - s_k^1| \leq \frac{1}{\eta} \sum_{x_j \in S}
 |w_j(h_k^1- \cos(P_{2,j})) \sref{x_j}| \nonumber\\
 & & \leq \frac{C}{\eta \, D}\sum_{x_j \in S}  d_k(\cdot,
 x_j)^3\,|\sref{x_j}| \leq O\left({\frac1D}\right),
 \end{eqnarray}
where $\eta$ is the constant provided by Theorem \ref{thm:auroux}.
Now let us check the $C^1$-norm of the error
\begin{eqnarray*}
 \left|dh_k^1 -d \left( \frac{s_k^1}{s_k^0} \right)\right| &=& \left|dh_k^1 -
 \frac{\nabla s_k^1 \cdot s_k^0 -s_k^1 \cdot \nabla s_k^0}{(s_k^0)^2} \right| =
 \left|\frac{1}{s_k^0} (dh_k^1 \cdot s_k^0 - \nabla s_k^1 +
 \frac{s_k^1}{s_k^0}\cdot \nabla s_k^0)\right| \\
 & \leq & \frac{1}{\eta} \Big(
 \sum |w_j|\,|dh_k^1-d \cos(P_{2,j})|\,|\sref{x_j}| + \\
 & + & \sum |w_j|\,
 |h_k^1-\cos(P_{2,j})|\,|\nabla \sref{x_j}| +
 \left| \frac{s_k^1}{s_k^0} - h_k^1 \right| |\nabla s_k^0| \Big)
 \leq \frac{1}{\eta} \frac{C}{D} =O\left({\frac1D}\right).
\end{eqnarray*}
Analogously, one can work out the $C^2$-norm of the error.
Choosing $D$ large, one can make this smaller than any $\gamma>0$.
\end{proof}

\section{Global transversality} \label{sec:global-trans}

In order to prove Theorem \ref{thm:main}, we use Proposition
\ref{prop:deform} to construct a sequence of maps $h_k:\SL \to
\R$. Then we use Proposition \ref{prop:nice_funct} to construct
two sequences of real asymptotically holomorphic sections $s_k^1$
and $s_k^2$ of $L^{\otimes k}$ such that
 $$
 \frac{\s_k^2}{\s_k^1}=
 \left(\frac{s_k^1+ i s_k^2}{s_k^1- i s_k^2}\right)\big|_\SL
 $$
is a Morse function isotopic to $(\cos (2h_k), \sin (2h_k))$. Here
$\s_k^1=s_k^1 - i s_k^2$ and $\s_k^2=s_k^1 +i s_k^2$. Let us check
the three transversality properties of Proposition
\ref{prop:get_tr} for this sequence in a tubular neighborhood of
fixed $g_k$-radius of $\SL$.

By the property given in the proof of Proposition
\ref{prop:nice_funct},
 $$
 \left|(h_k^1+i h_k^2) - \left( \frac{s_k^1 + i s_k^2}{s_k^0}\right)
 \big|_\SL \right|_ {C^2}
 $$
is small. Therefore $|\s_k^2|=|s_k^1+ i s_k^2| \geq |s_k^0|/2$
along $\SL$. In particular $\s_k^2$ is transverse to zero along
$\SL$. Analogously, $\s_k^1$ and $\s_k^1\oplus\s_k^2$ are
transverse to zero along $\SL$. Thus the first two properties of
Proposition \ref{prop:get_tr} hold in a tubular neighborhood of
$\SL$ of some uniform $g_k$-radius. Moreover, the sequence $\nabla
(h_k^1,h_k^2)$ is transverse to zero over $\SL$, so the sequence
$\nabla (\cos (2h_k), \sin (2h_k))$ is transverse to zero over
$\SL$. By Lemma \ref{lem:lag-trans}, this implies that $\nabla
(\s^2_k/ \s^1_k)$ is transverse to zero over a tubular
neighborhood of $\SL$ of fixed $g_k$-radius.

The transversality properties of Proposition \ref{prop:get_tr} are
$C^2$-open (see \cite{Au00} or Definition \ref{def:openn} below),
so we can follow Donaldson's ideas \cite{Do00} and perturb the two
sequences $\s_k^1,\s_k^2$ in order to achieve transversality over
all of $M$. However, the condition of $\SL$ projecting to a curve
is neither generic nor open (generically it projects to an open
set of $\CP^1$), so we cannot choose arbitrary Donaldson
perturbations. Hence some work has to be done. The method that we
follow is to perturb the quotient
 $$
 \frac{\s_k^2}{\s_k^1}= \frac{s_k^1+ i s_k^2}{s_k^1-is_k^2}
 $$
using only real perturbations of the sections $s_k^1$ and $s_k^2$.
Then the requirement that the sections be real forces the
image of the Lagrangian under the perturbed map to remain
in the unit circle.

An important remark to make at this point is that the property of being a
Lefschetz pencil is invariant under projective automorphisms of $\CP^1$,
and hence that the pair of sections $(s_k^1 - i s_k^2,s_k^1 + i s_k^2)$
defines a symplectic Lefschetz pencil if and only if the pair of real
sections $(s_k^1,s_k^2)$ defines a symplectic Lefschetz pencil. Therefore,
we only need to show how to perturb $s_k^1$ and $s_k^2$ among real sections
in order to ensure that they define a Lefschetz pencil.

We proceed as follows. Symplectically identify a neighboorhood $U$
of $\SL$ with a neighborhood $V$ of the zero section in $T^*\SL$.
We will define a number of structures over $T^*\SL$, understanding
that they translate to $U$ under this identification.
Let $e:T^*\SL \to T^*\SL$ be the involution of $T^*\SL$
defined by fiberwise
multiplication by $-1$, and denote by $\lambda$ the canonical
Liouville form on $T^*\SL$. We have $e^* \lambda = -\lambda$, and
therefore $e^* d\lambda = - d\lambda$. There is a special class of
compatible almost complex structures on $T^*\SL$ defined by
 $$
 \SJ(\SL) = \{ J\in \SJ(T^*\SL, d\lambda)\, : \, e^*J= -J \}.
 $$
It is simple to check that the set of all these almost complex
structures is contractible and non-empty. Let us call them
symmetric almost complex structures. Construct $g(u,v)=
d\lambda(u,Jv)$: then $e^* g=g$.

So far we have constructed a picture with a $\Z_2$-symmetry in
$T^*\SL$. Let us lift it to the bundle $L$. Over $T^*\SL$ we can
find a trivialization $L^{\otimes k}\cong T^*\SL \times \C$ in
which the connection is given by $\nabla = d+iA^0_k -ik\lambda$,
where $d+iA^0_k$ is a flat connection such that $e^*(A^0_k)=A^0_k$
and $|A^0_k|_{g_k}=O(k^{-1/2})$ everywhere. We can moreover assume
that the section $\s_{\SL,k}$ of $L^{\otimes k}$ over $\SL$
corresponds to the constant function $1$ in this trivialization.
Using this trivialization, the involution $e$ lifts to an
anti-complex isomorphism of line bundles
 \begin{eqnarray*}
 \hat{e}:\quad L^{\otimes k} &\to& L^{\otimes k} \\
      (x,z)&\mapsto& (e(x),\bar{z}).
 \end{eqnarray*}
The connection induced by $\hat{e}$ on $L^{\otimes k}$ is
$d-iA^0_k-ik\lambda$, which differs from $\nabla$ by
$O(k^{-1/2})$. Therefore, $\hat{e}^*(\bbd)=\bbd+O(k^{-1/2})$. We
say that a complex bundle $E$ is symmetric if $e$ admits a lift to
$E$ preserving its connection and complex structure. It is simple
to check that the tangent bundle and any bundle obtained from it
by tensor operations are symmetric bundles.

{}From all the previous remarks we get
\begin{lemma} \label{lem:symmetric}
Let $E$ be a symmetric bundle. Let $s_k$ be an asymptotically
holomorphic sequence of sections of the bundles $E\otimes
L^{\otimes k}$ with respect to a symmetric almost complex
structure. Then $\hat{e}^*(s_k)$ is also asymptotically
holomorphic. \newline Moreover, if $s_k$ is $\e$-transverse to
zero, then $\hat{e}^*(s_k)$ is also $\e$-transverse.
\end{lemma}

\begin{definition}\label{def:symmetric}
A section $s_k$ is called symmetric if $\hat{e}^*(s_k)=s_k$ in $U$ (and
no restriction outside $U$).
\end{definition}

Because we defined $\hat{e}$ using a trivialization of $L^{\otimes
k}$ over $\SL$ such that $\sigma_{\SL,k}=1$, we have

\begin{lemma}\label{lem:symmetric-real}
A symmetric section of $L^{\otimes k}$ is real.
\end{lemma}

Moreover, just by checking the construction in the proof of Lemma
\ref{lemma:local} we can show

\begin{lemma} \label{lem:symmetric-ref}
The reference sections constructed in Lemma \ref{lemma:local} can
be chosen to be symmetric whenever $x\in\SL$ if we fix a symmetric
almost complex structure. Moreover, the asymptotically
holomorphic sections constructed in Proposition
\ref{prop:nice_funct} can also be chosen symmetric.
\end{lemma}

\begin{proof}
The reference sections can easily be chosen to be symmetric.
Although a careful construction directly leads to symmetric
sections when $x\in\SL$, the most straightforward way is just to
choose reference sections
 $$
 \hat{\s}_{k,x}^{\mathrm{ref}}= \frac{\sref{x} + \hat{e}^*(\sref{x})}{2},
 $$
that still satisfy all the properties of the Lemma. In fact, we can also
use this method when $x\not\in\SL$, at the expense of obtaining sections
that are concentrated not only near $x$ but rather near the two points
$x$ and $e(x)$ (i.e., these sections satisfy suitably modified versions
of the second and third statements in Lemma \ref{lemma:local}).

The sections $s_k^1$ and $s_k^2$ constructed in Proposition
\ref{prop:nice_funct} are also symmetric. To check this, one just
has to check that the local sections $w_j \cos(P_{2,j})
\sref{x_j}$ are symmetric. To make sure that this condition is
fulfilled we just replace that section by
 $$
 \frac{w_j \cos(P_{2,j})\sref{x_j} + \hat{e}^*(w_j
 \cos(P_{2,j})\sref{x_j})}{2}
 $$
which is a section that coincides with the previous one over $\SL$
and satisfies the same properties (once again, since $x_j\in\SL$
and since $P_{2,j}$ are real polynomials, with a careful
construction this trick becomes unnecessary).
\end{proof}

{}From now on we make the following assumptions:

\begin{enumerate}
\item the almost complex structure $J$
is an extension to $M$ of a symmetric almost-complex structure in $U
\iso V \subset T^*\SL$.

\item
$J$ is integrable and standard in Darboux coordinates over
neighborhoods $B_{g_k}(p,c_0k^{1/2})$ of the critical points of $f$
(this assumption was already made in Section \ref{sec:morse}).

\item $\s_{\SL,k}$
is covariantly constant over a subset of $\SL$ containing
$B_{g_k}(p,c_0k^{1/2})$ for every critical point of $f$ (see the remark
after the statement of Lemma \ref{lem:trivialize}).
\end{enumerate}

Moreover, as in Lemma \ref{lem:symmetric-ref} we assume that we are dealing
with symmetric sections. In this context,
Remark \ref{rmk:holomorphic} still applies to our new symmetric
reference sections near the critical points of $f$, since under the
assumption that $\s_{\SL,k}$ is locally parallel the involution $\hat{e}$
is connection-preserving, so that averaging with $\hat{e}^*(\sref{x})$ does
not affect holomorphicity.

In Section \ref{sect:sym} we will prove the following result:

\begin{proposition} \label{prop:trans1}
Fix $\delta>0$ small enough. There exists an asymptotically
holomorphic sequence of symmetric sections $\tau_k^1\oplus
\tau_k^2$ such that $|\tau_k^1\oplus \tau_k^2|_{C^2}\leq \delta$
and satisfying
\begin{enumerate}
 \item $s_k^1+\tau_k^1$ is $\epsilon$-transverse to zero over $M$,
 for some uniform $\e>0$ (independent of $k$).
 \item $(s_k^1+\tau_k^1)\oplus(s_k^2+\tau_k^2)$ is
 $\e$-transverse to zero over $M$.
 \item $\bd\left(\dfrac{s_k^2+\tau_k^2}{s_k^1+\tau_k^1}\right)$ is
 $\e$-transverse to zero over the set $Z_{k,\e} =
 \{ p\in M\, : \, |s_k^1+\tau_k^1|\geq \e \}$.
\end{enumerate}
\end{proposition}

Applying this result to our sequence we get a new sequence that we
still denote by $s_k^1\oplus s_k^2$. The quotient of these two
sections $s_k^1$ and $s_k^2$, or equivalently the quotient
 $$
 \frac{\s_k^2}{\s_k^1}= \frac{s_k^1+ i s_k^2}{s_k^1-is_k^2},
 $$
satisfies the transversality estimates required of a symplectic
Lefschetz pencil. Moreover, because the perturbation
$\tau_k^1\oplus \tau_k^2$ can be chosen symmetric and $C^2$-small,
the restriction $(\s_k^2/\s_k^1)_{|\SL}$ remains a (circle-valued)
Morse function. Now we want to perturb the sections
$\s_k^1,\s_k^2$ in order to make the image of
$(\s^2_k/\s^1_k)|_\SL$ an embedded curve and obtain a Morse
function isotopic to the initial $h_k$.

Let $F_k=(\s^2_k/\s^1_k)|_\SL :\SL \to S^1$. Then $F_k$ is
$C^2$-close to $(\cos (2h_k), \sin (2h_k))$, therefore it can be
written in the form $(\cos (2\hat{h}_k),\sin (2\hat{h}_k))$, where the
function $\hat{h}_k$ is $C^2$-close to $h_k$. Extend $\hat{h}_k$
to the whole of $M$ keeping that
$|\hat{h}_k|_{C^2,g_k}=O(k^{1/2})$. Then define the pencil
  \begin{equation}\label{eq:embed}
  \tilde{F}_k=\frac{\s^2_k \cdot e^{\hat{h}_k/k}}{\s_k^1}
  \end{equation}
This is still asymptotically holomorphic since it is a
perturbation by $O(k^{-1/2})$ of the previous map $F_k$, and it still
satisfies the required transversality conditions.
Moreover, the restriction of $\tilde{F}_k$ to $\SL$ is equal to
 $$
 \exp(2i\hat{h}_k+ \hat{h}_k/k),
 $$
which is isotopic to $\exp(2ih_k+ h_k/k)$, i.e.\ the
composition of the initial $h_k$ with the embedding $\gamma (t)=
\exp(2it+t/k)$ of $[0,k^{1/2}]$ into $\C$. Also, observe that
these successive $C^2$-small perturbations affect in a controlled
manner the critical points of $\tilde{F}_k$ along $\SL$, which
remain within arbitrarily small $g_k$-distance of those of $h_k$.

We now just follow Donaldson's proof in \cite{Do00} to show that
$\tilde{F}_k$ defines a symplectic Lefschetz pencil. The proof has
three parts:
 \begin{itemize}
 \item We have to get the three transversality properties appearing in
   the statement of Proposition \ref{prop:get_tr}, either directly for $\tilde{F}_k$ or
   after composition with a projective automorphism of $\CP^1$.
 \item We have to check the local models near the base point
   set $\{p\in M: \s_k^1(p) \oplus \s_k^2(p)= 0 \}$.
 \item We have to check the local models close to the critical
   points of the map. By Donaldson's argument it automatically
   follows that the fibers are symplectic.
\end{itemize}

Proposition \ref{prop:trans1} covers the first part. The second
one follows Donaldson's argument without any adaptation. The third
one requires some work. We need the following auxiliary lemma:

\begin{lemma} \label{lem:quadratic}
 Let $p\in \SL$ be a critical point of $f$ (or equivalently, of
 $h_k$). Then we can arrange that $\s^2_k/\s^1_k$ is
 $J$-holomorphic over $B_{g_k}(p,c)$ for some fixed constant $c>0$, with a non-degenerate critical
point on $\SL$ within distance $c/2$ of $p$.
\end{lemma}

\begin{proof}
Let $x\in B_{g_k}(p, c_0k^{1/2})$. If $d_k(p,x)<k^{1/6}-c$, then Remark
\ref{rmk:holomorphic} and the assumptions made above imply that
$\sref{x}$ is holomorphic over $B_{g_k}(p,c)$.
If $d_k(p,x)>2k^{1/6}+c$ then the reference section $\sref{x}$
vanishes over $B_{g_k}(p,c)$. Finally, when $k^{1/6}-c
<d_k(p,x)<2k^{1/6}+c$, we can modify slightly the construction of the reference
sections in \eqref{eqn:sref}, using a cut-off function which equals $1$
inside a ball of $g_k$-radius $\frac12 k^{1/6}$ and $0$ outside of a ball
of $g_k$-radius $k^{1/6}-2c$ (instead of the usual distances $k^{1/6}$ and
$2k^{1/6}$); this again ensures the vanishing of $\sref{x}$ over $B_{g_k}(p,c)$.
If we use these reference sections (and linear combinations of them with
polynomial coefficients) throughout the construction (in particular in
Proposition \ref{prop:trans1}), then $\s_k^1$ and $\s_k^2$ are automatically
holomorphic over $B_{g_k}(p,c)$.

Recall that by adjusting the parameter $\delta$ in
Propositions \ref{prop:nice_funct} and
\ref{prop:trans1} we can assume that $(\s_k^2/\s_k^1)_{|\SL}$ is arbitrarily
$C^2$-close to $(\cos(2h_k),\sin(2h_k))$. In particular, because we have a
uniform estimate on the Hessian of $h_k$ at its critical points, by choosing
$\delta$ small enough in these two propositions we can ensure that the
critical points of $(\s_k^2/\s_k^1)_{|\SL}$ are non-degenerate and
lie within distance $c/2$ of those of $h_k$.
\end{proof}

With the above understood, if $p$ is a critical point of $h_k$
then $\s_k^2/\s_k^1=F_k(z)$ is holomorphic around $p$. So, when we
perturb $F_k$ in order to make the image of $\SL$ an embedded
curve, we can choose the function $\hat{h}_k$ in (\ref{eq:embed})
(with controlled $C^2$-norm, taking real values over $\SL$, and
such that $F_k=(\cos(2\hat{h}_k),\sin(2\hat{h}_k))$ over $\SL$) to
be holomorphic in a neighborhood of uniform $g_k$-radius of every
critical point. This means that throughout our construction we
have preserved holomorphicity near the critical points which lie
on $\SL$.

The last part of Donaldson's argument \cite{Do00} involves a perturbation
of the pencil map near its critical points in order to achieve the correct
(holomorphic quadratic) local model. However, we have already ensured
holomorphicity near the critical points in $\SL$ (any critical point
of $\tilde{F}_k$ along $\SL$ is also a critical point of $\hat{h}_k$ and
hence lies close to a critical point of $f$). Hence, this perturbation is
only needed at those critical points which lie outside of $\SL$. Because
$\nabla h_k$ is $\eta$-transverse to zero over $\SL$ we can
find a uniform lower bound on the $g_k$-distance between such a critical
point and $\SL$. Therefore the required perturbation can be carried out without
affecting the behavior of $\tilde{F}_k$ along $\SL$, thus completing the
proof of Theorem
\ref{thm:main}.

\section{Globalization preserving symmetry.} \label{sect:sym}

\subsection{Proof of the first part of Proposition \ref{prop:trans1}.}
\label{sub:first} We have to follow Donaldson's ideas in
\cite{Do96}. Let us review how to make a section $s_k$ transverse
to zero. The steps to follow are

\begin{enumerate}
\item Given a point $x\in M$, we use the chart $\Psi_{k,x}$
 provided by Lemma \ref{lem:darboux} to trivialize in an
 approximately holomorphic way a neighborhood of fixed $g_k$-radius
 of $x$. We also trivialize the bundle $L^{\otimes k}$ in that
 neighborhood using the section $\sref{x}$ given by Lemma
 \ref{lemma:local}. This allows us to construct an asymptotically
 holomorphic function
   \begin{eqnarray*}
    f_{k,x}:B(0,1)\subset \C^n  & \to & \C \\
    y & \mapsto &
    \frac{s_{k}(\Psi_{k,x}^{-1}(y))}{\sref{x}(\Psi_{k,x}^{-1}(y))}\, .
   \end{eqnarray*}
\item We perturb the function to a new one which is transverse to
 zero by adding a small constant $w$ (this amounts to adding
 $w\sref{x}$ to $s_k$). The precise result is

\begin{lemma}[Donaldson] \label{lem:trans_real}
There exists an integer $p$ depending only on the dimension $n$,
with the following property: let $\delta$ be a constant with
$0<\delta<1/2$ and $\s= \delta (\log(\delta^{-1}))^{-p}$. Let
$f:B(0,11/10) \subset \C^n \to \C$ satisfy
  $$
  |f| \leq 1, ~~ |\bbd f|\leq \s, ~~ |\nabla \bbd f| \leq \s.
  $$
Then there exists a complex number $w$ with $|w|\leq \delta$, such
that $f-w$ is $\sigma$-transverse to zero in the unit ball of
$\C^n$.
\end{lemma}

\item The third part of the argument consists of a method to achieve
transversality in the entire manifold from the local
transversality argument. We start by fixing a covering of the
manifold by a finite number of Darboux charts $\{U_i\}_{i=1}^N$.
It is enough to know how to get transversality in one of them. For
we can start by perturbing $s_k$ to some $s_k^1=s_k+\sigma_k^1$
which is $\e_1$-transverse to zero in $U_1$ and then we perturb
every $s_k^i$ to $s_k^{i+1}=s_k^i+\sigma_k^{i+1}$ which is
$\eta_i$-transverse to zero in $U_{i+1}$. If $s_k^i$ is
$\e_i$-transverse to zero in $\cup_{j=1}^i U_j$ and the
perturbation is small enough then $s_k^{i+1}$ is still
$\e_i/2$-transverse to zero in $\cup_{j=1}^i U_j$. This is
because transversality is a $C^1$-open
property, in the following sense:

\begin{definition} \label{def:openn}
A property $\SP(\e,x)$ of a sequence of sections $s_k:M\to E_k$ is
$C^r$-open if there exists a uniform constant $c_0>0$ such that if
$s_k$ has the property $\SP(\e,x)$ then any $s_k'$ such that
$|s_k(x)-s_k'(x)|_{C^r}\leq c_0 \e$ has the property
$\SP(\e/2,x)$.
\end{definition}

Choosing $\e_{i+1}=\min\{\e_i/2,\eta_i\}$, $s_k^{i+1}$ is
$\e_{i+1}$-transverse to zero in $\cup_{j=1}^{i+1} U_j$.

The difficult part is to get transversality in a single Darboux
chart. Take the lattice $(c\Z)^{2n}\subset \R^{2n}$, for some
small $c>0$. The image of the lattice in $M$ will be denoted by
$\Gamma$. We perturb our section $s_k$ with sections centered at
the different points of the lattice $\Gamma$. An argument of Donaldson
\cite{Do96} shows a way of doing this which yields uniform
transversality over the entire Darboux chart, with a
transversality amount independent of $k$, in spite of the
increasing number of lattice points needed to cover the rescaled
Darboux chart. The crucial ingredients here are: (a)
transversality is an open property, and (b) the sections
$\sref{x}$ decay exponentially away from $x$. These two features
make it possible to keep track of the mutual interference created
by different perturbations. We do not provide the details, which
can be found in \cite{Do96}.
\end{enumerate}

In our case, we have to adapt the argument to perturb an
asymptotically holomorphic sequence of {\em symmetric\/} sections
$s_k^1$ to make it transverse to zero. This means that we want the
construction to be invariant by the involution $\hat{e}$.

As above, we start by covering the manifold by a finite number of
Darboux balls $U_j$, chosen in such a way that the Darboux balls
which intersect the Lagrangian submanifold $\SL$ satisfy
$e(U_j)=U_j$. As before, we need to get transversality for a
single chart. If $U_j$ does not intersect $\SL$ then we use
Donaldson's argument without modification. If the chart intersects
$\SL$ then we have to be more careful. We will get transversality
in pairs of small balls $V$ and $e(V)$ at the same time. Let us
review the adaptations needed in each of the three parts of the
argument:
 \begin{enumerate}
  \item For the first part, there is no change. We trivialize our section
  $s_k^1$ in a neighborhood of a point $x$, using the trivializing section
  $\sref{x}$ to obtain a complex-valued function $f$ defined by
  $$ f= \frac{s_k^1}{\sref{x}}.$$
\item We also follow the structure of the second part.
Now, the section $w \sref{x}+ \bar{w} \hat{e}^*{\sref{x}}$ is symmetric
for any complex number $w$, and can therefore be used to locally
perturb $s_k^1$. After division by $\sref{x}$, replacing
$s_k^1$ by $s_k^1-w \sref{x}-\bar{w} \hat{e}^*{\sref{x}}$ corresponds to a
perturbation of $f$ to the new function
$$ f_{(w)}= f-w- \bar{w}\frac{\hat{e}^*{\sref{x}}}{\sref{x}}, $$
Now we apply the following
\begin{proposition} \label{local-trans-equiv}
 Let $C>0$, $\kappa>0$.
 There exists an integer $p$ and a fixed $\delta_0>0$
 (depending only on $C$, $\kappa$ and the dimension)
 verifying the following property:
 for $0<\delta<\delta_0$, let $\sigma=\delta(\log (\delta^{-1}))^{-p}$.
 Let $f$, $h$ be two complex-valued functions defined over the ball $B^+=
 B(0,\frac{11}{10}) \subset \C^n$ satisfying the following bounds over $B^+$,
\begin{eqnarray*}
  &  |f|\leq 1,\quad |\bbd f|\leq \s,\quad |\nabla \bbd f|\leq \s,  &\\
  &  |h|\leq 1-\kappa,\quad|d h|\leq C, \quad
  |\bbd h|\leq \s,\quad |\nabla \bbd h|\leq \s. \\
\end{eqnarray*}
 Then there exists $w \in \C$ with $|w|<\d$ such that $f-w-\bar{w}\, h$ is
 $\s$-transverse to $0$ over the unit ball in $\C^n$. The same result holds
 for one-parameter families of
 functions $f_t$, $h_t$ depending continuously on $t\in [0,1]$, where
 we obtain a continuous path $w_t \in \C$, with $|w_t|<\d$.
\end{proposition}

The proof of this Proposition will be given in Section \ref{sec:parametric};
it is essentially a careful modification of Donaldson's local Lemma.

Recall from the beginning of Section \ref{sec:global-trans} that the desired
transversality property already holds over a neighborhood of $\SL$ of
uniform $g_k$-radius, say $\rho>0$, and will automatically be preserved by all
perturbations of $s_k^1$ that are sufficiently $C^1$-small. Therefore,
in the construction
we only need to consider the case where $x$ is at distance at least $\rho$
from $\SL$, and rescaling the coordinates by a constant factor if necessary
we can assume that the ball $B^+$ does not intersect the
$(\rho/2)$-neighborhood of $\SL$. This implies that
$h=\hat{e}^*(\sref{x})/\sref{x}$ satisfies a bound of the form $|h|\le
1-\kappa$ over $B^+$ for a uniform constant $\kappa>0$, thus allowing us
to apply Proposition \ref{local-trans-equiv}.

Then the new section $s_k^1 - w\sref{x}- \bar{w}\hat{e}^*(\sref{x})$ is
$\s$-transverse to zero in a ball $B_x$ of fixed radius around
$x$. Moreover, it is symmetric as desired, and therefore by symmetry it is
also $\s$-transverse to zero over $e(B_x)$.

\item The third part remains unchanged except for the fact
that we choose the lattice $\Gamma$ to be $e$-invariant, and the
transversality is obtained simultaneously over pairs of balls centered
at points $x$, $e(x)$, for all $x\in\Gamma$ which lie at $g_k$-distance
at least $\rho$ from $\SL$ (while transversality already holds over the
$\rho$-neighborhood of $\SL$ by assumption).
\end{enumerate}

\subsection{Proof of the second part of Proposition \ref{prop:trans1}.}
\label{sub:second}

Now we have to get transversality for the sequence $s_k^1\oplus
s_k^2$ while preserving $e$-invariance. Instead of using the
generalization of the local transversality Lemma
\ref{lem:trans_real} to maps from $\C^n$ to $\C^r$ with $r>1$
obtained by Donaldson in \cite{Do00}, we shall adapt the ideas of
\cite{Au97} to our $e$-invariant setting. This is because we
haven't developed the extension of Proposition
\ref{local-trans-equiv} for maps from $\C^n$ to $\C^r$. Such a
result should hold, but it is easier for the purposes of this
article just to carefully use the case $r=1$.

Because of the transversality of the sequence $s_k^1$ achieved by means
of a symmetric perturbation in Subsection \ref{sub:first},
the zero set $Z(s_k^1)$ of the section $s_k^1$ is a
symplectic submanifold for $k$ large enough, and the tangent
spaces to $Z(s_k^1)$ are within $O(k^{-1/2})$ of being complex
subspaces of $(TM,J)$. We perturb $s_k^2$ in order to make its
restriction to the submanifold $Z(s_k^1)$ transverse to zero. The
argument is similar to the one used to get transversality of
$s_k^1$ all over $M$, however some details concerning the fact
that the submanifold $Z(s_k^1)$ changes with $k$ have to be taken
into account. We refer the reader to \cite{Au97} for more details.
In our case, the only required adaptation is to impose symmetry
for all the perturbations. However, this works exactly in the same
way as in Subsection \ref{sub:first}, and only requires us to use
Proposition \ref{local-trans-equiv}.

Once we have perturbed the sequence $s_k^2$ into one that is
$\e$-transverse to zero on $Z(s_k^1)$, we conclude with the
following linear algebra lemma:

\begin{lemma}[Subsection 3.6 in \cite{Au97}] \label{lemma:adding_tr}
Let $s_k^1\oplus s_k^2$ be an asymptotically holomorphic sequence
of sections of the bundles $E_k=(E_1 \oplus \C) \otimes L^{\otimes
k}$. Assume that $s_k^1$ is $\e$-transverse to zero on an open set
$U$. Assume that $s_k^2$ is $\e$-transverse to zero on
$Z(s_k^1)\cap U$. Then $s_k^1\oplus s_k^2$ is $\e'$-transverse to
zero over $U$, where $\e'>0$ depends only on $\e$ and the uniform
bound on $|\nabla\nabla(s_k^1\oplus s_k^2)|$, but not on $k$.
\end{lemma}

\subsection{Proof of the third part of Proposition \ref{prop:trans1}.}
\label{sub:third} {}From the previous steps of the argument we may
assume that there exists an uniform $\e>0$ for which
 \begin{itemize}
  \item $s_k^1\oplus s_k^2$ is $\e$-transverse to zero over $M$.
  \item $s_k^1$ is $\e$-transverse to zero over $M$. \item $s_k^2$
  is $\e$-transverse to zero over $M$.
 \end{itemize}
So, our only remaining task is to modify $s_k^1\oplus s_k^2$ by
symmetric perturbations in order to achieve transversality of the
sequence $\bd(s_k^2/s_k^1)$ over the set $Z_{k,\e}$. As in
Subsection \ref{sub:first}, the $C^2$-openness of this condition
(see Definition \ref{def:openn}) makes it possible to reduce the
problem to a fixed Darboux chart $U$; if $U$ does not intersect
the Lagrangian submanifold $\SL$, then we just refer to the
standard proof in \cite{Do00}. In the other case, let $U$ be a
Darboux chart intersecting $\SL$, such that $e(U)=U$, and fix a
set of sections $\{ \alpha_j \}_{j=1}^n$ of the topologically
trivial bundle $(T^*U)^{1,0}$, defining an orthonormal basis at
every point of $U$. Note that we do not impose any integrability
condition on the sections.
Define the line subbundles $A_r= \C\alpha_r\subset
(T^*U)^{1,0}$, and the corresponding filtration $E_j=
\bigoplus_{r=1}^j A_r$. The decomposition
$(T^*U)^{1,0}=E_n=\bigoplus_{r=1}^n A_r$ allows us to define
projection operators $\pi_{r}: E_n \to A_r$ and $\pi_{\le r}: E_n
\to E_r$. We consider the section $\bd(s_k^2/s_k^1)|_{U}:U \to
E_n$, and introduce the following definition:

\begin{definition} \label{def:sharp}
A section $s_k^1\oplus s_k^2$ of $\C^2\otimes L^{\otimes k}$ is
$(\e,E_r)$-sharp at a point $x\in U$ if the section $\pi_{\le
r}(\partial(s_k^2/s_k^1))$ of $E_r$ is $\e$-transverse to zero at
$x$.
\end{definition}

Note that this property is $C^2$-open as defined in Definition
\ref{def:openn}. This allows us to use the three step process
described in Subsection \ref{sub:first} to achieve sharpness at every
point of $U\cap Z_{k,\epsilon}$.
Moreover, $E_n$-sharpness is equivalent to transversality to zero
for $\partial(s_k^2/s_k^1)$ and $E_0$-sharpness is a void
condition.

The scheme of the proof is to assume $(\e,E_r)$-sharpness and
then manage to achieve $(\e',E_{r+1})$-sharpness.
Denote by $W_{k,r}$ the zero set of $\pi_{\le r}(\partial
(s_k^2/s_k^1))$, which is a symplectic submanifold for $k$ large
enough. We proceed by making the following definition:

\begin{definition} \label{def:Asharp}
A section $s_k^1\oplus s_k^2$ is $(\e,A_{r+1})$-sharp at a point
$x\in U\cap W_{k,r}$ if the restriction to $W_{k,r}$ of the
section $\pi_{r+1}(\partial(s_k^2/s_k^1))$ of $A_{r+1}$ is
$\e$-transverse to zero at $x$.
\end{definition}

This is again a $C^2$-open property under the assumption of
$E_{r}$-sharpness (because $W_{k,r}$ then depends continuously on
$s_k^1\oplus s_k^2$). By Lemma \ref{lemma:adding_tr}, we know that
$E_{r+1}$-sharpness on $U$ follows from $A_{r+1}$-sharpness on
$W_{k,r}\bigcap U$.

We analyze the three steps needed for the argument:
\begin{enumerate}
\item Fix a point $x\in W_{k,r}$. We use Lemma \ref{lem:darboux}
to define a Darboux chart $\Psi_{k,x}$ centered at $x$.
After applying a linear transformation in $U(n)$, we can ensure
that, at the point~$x$,
 $$
 \alpha_j(x)= (\Psi_{k,x}^* (\bd z_j))(x).
 $$
Also, recall from the beginning of Section \ref{sec:global-trans}
that we can assume $x$ lies at a $g_k$-distance bigger than a
uniform constant from $\SL$ (indeed, the desired tranversality
property over a small tubular neighborhood of $\SL$ follows from
the estimates achieved in Section \ref{sec:trans_1}). This implies
the existence of a uniform lower bound on at least one of the
coordinates at $e(x)$, say $z_{j_0}$. Let
$\gamma_x=\Psi_{k,x}^*(z_{j_0})(e(x))$, and for $w\in\C$ define
 $$
 s^2_{k;(w)}=s_k^2+w\Psi_{k,x}^*(z_{r+1}(1-\gamma_x^{-1}z_{j_0})^2)\sref{x}.
 $$
Since we have a uniform upper bound on $\gamma_x^{-1}$, the perturbations
we add to $s^2_k$ are asymptotically holomorphic and have uniform exponential
decay away from $x$, as required for the general construction.
For small enough values of $w$ the submanifold
$W_{k,r;(w)}=\{\pi_{\le r}(\partial(s^2_{k;(w)}/s^1_k))=0\}$ is a small
deformation of $W_{k,r}$, and can be viewed as the image by the exponential
map of the graph of a section of the normal bundle to $W_{k,r}$. With this
understood, we can define
 $$
 \zeta_{k,x}=\frac{\partial}{\partial w}_{|w=0}
 \bigl(\pi_{r+1}(\partial(s^2_{k;(w)}/s^1_k))|_{W_{k,r;(w)}}\bigr).
 $$
Observe that $\zeta_{k,x}(x)=(\alpha_{r+1}(x)\sref{x}(x)/s^1_k(x))$ is uniformly
bounded from below, and so there is a uniform radius $\rho>0$ such that
$\zeta_{k,x}$ is bounded from below by a uniform constant over $W_{k,r}\cap
B_{g_k}(x,\rho)$. Moreover, by construction the perturbation we add to
$s_k^2$ vanishes up to order $2$ at $e(x)$, so $\zeta_{k,x}(e(x))=0$.
Hence, decreasing $\rho$ if necessary (but preserving uniformity) we can
assume that a bound of the form $|\hat{e}^*\zeta_{k,x}|\le (1-\kappa)
|\zeta_{k,x}|$ holds at all points of $W_{k,r}\cap B_{g_k}(x,\rho)$.

Finally, recall from \cite{Au97} that, decreasing $\rho$ if necessary, we
can trivialize the submanifold $W_{k,r}$ in the ball of radius $2\rho$
around $x$. More precisely, using Lemma~4 of \cite{Au97} we can get an
approximately holomorphic chart $\Phi_{k,x} : B(0,\frac{11}{10})
\subset \C^{n-r} \to W_{k,r}$. Using these trivializations we can introduce
a function $f_{k,r+1}:B(0,\frac{11}{10}) \subset \C^{n-r} \to \C$, defined by
 $$
 f_{k,r+1}=\Phi_{k,x}^*\left(\frac{\pi_{r+1}(
 \partial( s^2_k/s^1_k))}{\zeta_{k,x}}\right),
 $$
describing the restriction to $W_{k,r}$ of the component of
$\partial(s_k^2/s_k^1)$ along $\alpha_{r+1}$, and satisfying
asymptotically holomorphic estimates.

\item If $x$ lies on $\SL$ or within a small tubular neighborhood of uniform
$g_k$-radius, then the transversality estimates obtained on $\SL$ at the
beginning of Section \ref{sec:global-trans} imply the desired result, and
no work is required. Otherwise,
we apply Proposition \ref{local-trans-equiv} to get a constant
$w$ such that
 $$
 f_{k,r+1;(w)}=f_{k,r+1}-w- \bar{w}\,\Phi_{k,x}^*\left(
 \frac{\hat{e}^*(\zeta_{k,x})}{\zeta_{k,x}}\right)
 $$
is $\sigma$-transverse to zero over the
unit ball. Now we define
 \begin{equation}\label{eqn:eqn}
 \s_k^2= w\Psi_{k,x}^*(z_{r+1}(1-\gamma_x^{-1}z_{j_0})^2) \sref{x} +
 \hat{e}^*\Bigl(w\Psi_{k,x}^*(z_{r+1}(1-\gamma_x^{-1}z_{j_0})^2)\sref{x}\Bigr).
 \end{equation}
Clearly, the new sequence $s_k^1\oplus (s_k^2 - \s_k^2)$
corresponds to the desired perturbation (i.e., replacing
$f_{k,r+1}$ with $f_{k,r+1;(w)}$ in the trivialization), up to
higher order terms (bounded by $O(w^2)$), which are negligible,
i.e., much smaller than $\sigma$ (recall that
$\s=\d(\log(\d^{-1}))^{-p}$ is the constant appearing in
Proposition \ref{local-trans-equiv}). This gives us
$(\s,A_{r+1})$-sharpness over a ball of fixed $g_k$-radius in
$W_{k,r}$. Also, by construction the perturbation (\ref{eqn:eqn})
is $\hat{e}$-invariant.

\item The final step (patching together the various local perturbations
to achieve sharpness globally) is the same as before: we again choose the
lattice $\Gamma$ to be $e$-invariant, and in order to be able to apply
Proposition \ref{local-trans-equiv} we discard those points of
$\Gamma$ which lie inside the small neighborhood of $\SL$ over which
the initially given sections already satisfy all desired properties.
\end{enumerate}

\section{The parametric case} \label{sec:parametric}

Theorem \ref{thm:main_par} is proved using the same strategy as
Theorem \ref{thm:main}, replacing all the objects in the above
argument by families indexed by the parameter $t\in [0,1]$. The
existence of families of Darboux charts $\Psi_{k,x,t}$ and reference
sections $\sref{x,t}$ that depend continuously on $t$, as well as
one-parameter versions of the various standard results in
approximately holomorphic geometry, are well-known \cite{Au97}
\cite{Do00}. It is worth noting that, since the isotopy
$(\SL_t)_{t\in [0,1]}$ is Hamiltonian and hence induced by a family
of global symplectomorphisms $\psi_t$ of $M$, after pullback via
$\psi_t$ we can assume that the Lagrangian submanifold is fixed
(while all the other auxiliary data in the construction depend on
$t$). Although the standard results can be used to achieve
transversality along submanifolds that depend on $t$ (see e.g.\ the
arguments in \cite{Au97}), keeping $\SL$ fixed makes it easier to
adapt the argument to the parametric case.

In any case, all that is required is to check that, at every stage of the construction
of adapted pencils (throughout Sections \ref{sec:morse}--\ref{sect:sym}),
it is possible to interpolate between the choices of sections and
perturbations which lead to the given pencils $\phi_{k,0}$, $\phi_{k,1}$ via
one-parameter families of objects that depend continuously on $t$.
The two places in the argument where the adaptation to the parametric
case is not completely straightforward are the existence of the real
trivialization $\s_{\SL,k}$ (Lemma \ref{lem:trivialize}), and the
possibility of achieving estimated transversality via symmetric perturbations
(Proposition \ref{local-trans-equiv}).

The existence of a continuous family of sections $\s_{\SL_t,k}$ of
$L^{\otimes k}$ over $\SL_t$ which satisfy the properties of Lemma
\ref{lem:trivialize} and extend the given choices for $t=0$ and
$t=1$ is guaranteed by the requirement that $\pi_1(\SL_t)=1$.
 This ensures that the non-vanishing
sections $\s_{\SL_0,k}$ and $\s_{\SL_1,k}$ used to build the pencils
$\phi_{k,0}$ and $\phi_{k,1}$ belong to the same homotopy class.
Even if $\pi_1(\SL_t)\neq 1$, Theorem \ref{thm:main_par} still holds when
the isotopy $\{\SL_t\}_{0\le t\le 1}$ is generated by a family of global
symplectomorphisms of $M$, provided that the linking numbers of the fibers
of the pencil with all closed loops on the Lagrangian submanifold are
assumed to be equal for $\phi_{k,0}$ and $\phi_{k,1}$ (or equivalently,
assuming that $\s_{\SL_0,k}$ and $\s_{\SL_1,k}$ are homotopic). Note that, if the
isotopy $\{\SL_t\}$ is not generated by global symplectomorphisms, then
Stokes' theorem precludes the
existence of a continuous family of non-vanishing sections $\s_{\SL_t,k}$
such that $|\nabla\s_{\SL_t,k}|=O(k^{-1/2})$.

Finally, concerning Proposition \ref{local-trans-equiv}, we give below a
proof that adapts easily to the parametric case (in the non-parametric
case it would be easier to follow the argument in \cite{Do96} and
show that $w$ can be chosen real).

\subsection*{Proof of Proposition \ref{local-trans-equiv}.}
\begin{proof}
  As in \cite{Do96,Au97}, we can approximate $f$ and $h$ by
  polynomials $p$ and $q$ such that
\begin{eqnarray} \label{eqn:above}
 && |f-p|_{C^1}\leq c\,\s,~~~~ \deg(p) \leq  c\,\log(\s^{-1}), \nonumber \\
 && |h-q|_{C^1}\leq c\,\s,~~~~ \deg(q) \leq  c\,\log(\s^{-1}),
\end{eqnarray}
in the unit ball, for $c>0$ a constant. Now define the function
 $$
 s(z,w)= p(z) - w -\bar{w}\, q(z).
 $$
For $\s$ small enough (choosing $\d_0$ appropriately), we get that
$|q|\leq 1-\frac{\kappa}{2}$. Therefore
 $$
 \frac{\partial s}{\partial w}= -\id - \frac{\partial \bar{w}}{\partial w}\, \cdot q
 $$
is an invertible matrix (here $\partial s/\partial w$ means the Jacobian
matrix of the function $s$ with respect to the two variables
$\mathrm{Re}(w),\,\mathrm{Im}(w)$). Its inverse satisfies
 \begin{equation}\label{eqn:111}
  \left|\left(\frac{\partial s}{\partial w}\right)^{-1}\right| \leq
  2\kappa^{-1}.
 \end{equation}
We can apply the implicit function theorem to $s(z,w)$ to obtain a
(unique) function $w=w(z)$ such that $s(z,w(z))= 0$. This function
is defined over $B$ and, taking derivatives, we get
 $$
  \frac{\partial s}{\partial z} +  \frac{\partial s}{\partial w}
  \frac{dw}{dz}=0.
 $$
So
 $$
  \frac{dw}{dz} = - \left(\frac{\partial s}{\partial w}\right)^{-1} \,
  \frac{\partial s}{\partial z}
 $$
along the graph of $w=w(z)$. This implies that
 \begin{equation}\label{eqn:222}
  \left|\frac{dw}{dz}\right| \leq 2\kappa^{-1} \bigl|l(z)\bigr|,\quad
  \mathrm{where} \quad l(z)=\frac{\partial s}{\partial z} \, (z,w(z)).
 \end{equation}
Recall that $\partial s/\partial z$ is the Jacobian matrix of the partial
derivatives of $s$ with respect to $\mathrm{Re}(z)$ and $\mathrm{Im}(z)$,
which we evaluate along the graph of $w=w(z)$. Now consider the set
 $$
 Y_{l,\s}= \{ z\in B: \, |l(z)|\le \s\}.
 $$
Given constants $c_1,c_2>0$,
denote by $N_{c_2\s}(w(Y_{l,c_1 \s}))$ the $c_2\s$-neighborhood
of the image $w(Y_{l,c_1\s})\subset \C$. Suppose that $w_0$ is a
point outside of this neighborhood. Let us see that $s(z,w_0)$ is
$c_0\,\s$-transverse to zero over $B$, for some constant $c_0>0$.

Take a point $z\in B$ with $|s(z,w_0)-s(z,w(z))|=|s(z,w_0)|<
\frac{1}{2}\kappa c_2\,\s$. Then $|w(z)
-w_0| < c_2\,\s$ by \eqref{eqn:111}. Therefore $z$ is not in
$Y_{l,c_1 \s}$ and so $|l(z)|\geq c_1\s$. This is rewritten as
 $$
 \left|\frac{\partial s}{\partial z}(z,w(z))\right| \geq c_1\s.
 $$
At this point we apply the bound
 $$
 \left|\frac{\partial ^2s}{\partial z\,\partial w}\right|=
 \left|\frac{dq}{dz}\right| \leq C+1
 $$
for $\s$ small enough. Decreasing $c_2$ if necessary we can assume that $c_1>2(C+1)c_2$,
and then
 $$
 \left|\frac{\partial s}{\partial z}(z,w_0)\right| \geq c_1\s - (C+1)c_2\s \geq
 \frac{c_1}2 \s.
 $$
 Therefore $s(z,w_0)$ is
$c_0\,\s$-transverse to zero with $c_0=\min\{\frac{1}{2}c_1,\frac{1}{2}\kappa c_2\}$.

Finally, recall that transversality is a $C^1$-open property, and that
(assuming $|w_0|<1$) the bounds \eqref{eqn:above} imply that $|s(z,w_0)-
(f-w_0-\bar{w}_0\, h)|_{C^1}\leq 2c\,\s$.
Since  $s(z,w_0)$ is $c_0\s$-transverse to zero, it follows that
$f-w_0-\bar{w}_0\, h$ is $c_0\s/2$-transverse to zero if the constant
$c_0$ is large enough compared to $c$.

It remains to see that there are points in $B(0,\d)-
N_{C\s}(w(Y_{l,C\s}))$ for $C$ a large constant. For this we use
the following result (Proposition 25 of \cite{Do96}):

\begin{lemma} \label{lem:Y}
 Let $F:\R^m \to \R$ be a polynomial function of degree $d$, and
 let $S(\theta)=\{x\in\R^m :  \left|x\right| \leq 1,\ F(x) \leq
 1+\theta \}$. Then for arbitrarily small $\theta >0$ there exist
 constants $C$ and $\nu$ depending only on the dimension $m$ such
 that $S(0)$ may be decomposed into pieces $S(0)=S_1 \cup \cdots
 \cup S_A$, where $A\leq C d^\nu$, in such a way that any pair of
 points in the same piece $S_r$ can be joined by a path in
 $S(\theta)$ of length less than $Cd^\nu$.
\end{lemma}

Now it is necessary to give a more explicit description of the
function $w(z)$. Since it solves the equation
$p(z)=w(z)+\bar{w}(z)\, q(z)$, we have
 $$
  w= \frac{p-\bar{p}q}{1-|q|^2}.
 $$

Therefore $w(z) (1-|q|^2)$ is a polynomial function of the $2n$ real
variables $\mathrm{Re}(z_j),\mathrm{Im}(z_j)$. Hence the condition
$|l(z)|\leq C\s$ can be rewritten as
$|l(z)|^2(1-|q|^2)^2\leq C^2\s^2(1-|q|^2)^2$,
which is equivalent to
 \begin{equation}\label{eq:lbound}
  \left| (1-|q|^2)\,\frac{dp}{dz} - (1-|q|^2)\,\bar{w}(z)\,
   \frac{dq}{dz} \right|^2 - C^2 \s^2 (1-|q|^2)^2+1 \leq 1
 \end{equation}
Let $F$ be the polynomial on the left hand side of \eqref{eq:lbound}, so that $S(0)$
is $Y_{l,C\s}$. By Lemma \ref{lem:Y}, there is a polynomial $P(d)$
such that $S(0)$ can be decomposed into at most $P(d)$ subsets,
such that two points in the same subset can be joined by a path of
length at most $P(d)$ inside $Y_{l,2C\s}$ (choose
$\theta=3C^2\kappa^2\s^2/4 \leq 3C^2\s^2 (1-|q|^2)^2$, since $(1-|q|^2)^2>\kappa^2/4$).
Using the bound \eqref{eqn:222} on $|dw/dz|$, this implies that the image
$N_{C\s}(w(Y_{l,C\s}))$ is contained in the union of $P(d)$ discs
of radius at most $(4\kappa^{-1}P(d)+1) C\s$
(see the argument in \cite[pages 976-977]{Au97}).
This proves that for some large fixed $p>0$, if $\s=\d(\log
(\d^{-1}))^{-p}$, then there is a connected component of $B(0,\d)
-N_{C\s}(w(Y_{l,C\s}))$ of area more than $0.9 \pi\d^2$. This is
enough to finish the proof. The argument also applies in the parametric
case, as shown in \cite{Au97}.
\end{proof}

\section{From matching paths to Lagrangian submanifolds}
\label{sec:matching}

Theorem \ref{thm:main} shows how to place a Lagrangian submanifold
in special position with respect to a Lefschetz pencil. In the converse
direction, it is natural to ask when it is possible to reconstruct a
Lagrangian submanifold from a given Lefschetz pencil and a given path in
the base $\CP^1$. We will restrict the discussion to the case of
Lagrangian spheres, following arguments due to Donaldson and Seidel. Many
of the ideas can be adapted for general Lagrangian submanifolds,
but some technicalities are involved. According to Theorem
\ref{thm:main} a Lagrangian sphere can be realized as the lift of a
path joining two critical points.

Recall that, outside of the base points and critical points, a symplectic
Lefschetz pencil carries a natural horizontal distribution given by the
symplectic orthogonal to the fiber. Parallel transport along an arc
$\gamma:[0,1]\to \CP^1$ which avoids the critical values induces a
symplectomorphism between the smooth fibers above the endpoints. However,
if one of the endpoints $\gamma(0)$ is a critical value, then one can
associate to $\gamma$ a Lagrangian ``vanishing disc'' $D\subset M$ (also
called ``Lefschetz thimble''), which is the set of all the points in the
fibers above $\gamma$ for which parallel transport converges to the
critical point in the fiber above $\gamma(0)$. The boundary $S=\partial D$
is a Lagrangian sphere in the fiber above $\gamma(1)$, called the
{\it vanishing cycle} associated to $\gamma$.
With this understood, we give the following

\begin{definition} \label{def:matching}
A matching path for a symplectic pencil $\phi$ is an embedded
curve $\gamma: [0,1] \to \CP^1$ such that:
\begin{enumerate}
 \item the only critical values of the pencil which lie on the curve
  $\gamma$ are the endpoints $\gamma(0)$ and $\gamma(1)$;
 \item the vanishing cycles $S_0, S_1\subset \phi^{-1}(\gamma(\frac12))$
  associated with the arcs $\gamma|_{[0,1/2]}$ and $\gamma|_{[1/2,1]}$
  are Hamiltonian isotopic inside the fiber $\phi^{-1}(\gamma(\frac12))$
  through a Hamiltonian isotopy with support away from the base point set.
  \end{enumerate}
\end{definition}

A matching path can be used to construct a smoothly embedded
$n$-sphere $S_\gamma\subset M$ in the following manner. Let
$(\Psi_s)_{s\in[0,1]}$ be a Hamiltonian isotopy in the fiber
$\phi^{-1}(\gamma(\frac12))$ connecting $S_0$ and $S_1$ (so
$\Psi_0=\Id$ and $\Psi_1(S_0)=S_1$), and let $\chi:[0,1]\to[0,1]$
be a smooth function such that $\chi|_{[0,1/3]}=0$ and
$\chi|_{[2/3,1]}=1$. For $t\in[0,1]$, denote by
$\SP_t:\phi^{-1}(\gamma(\frac12))\to\phi^{-1}(\gamma(t))$ the map
induced by parallel transport along $\gamma$. Then we let
$S_\gamma=\bigcup_{t\in[0,1]}\SP_t(\Psi_{\chi(s)}(S_0))$.

Although the sphere $S_\gamma$ obtained in this way coincides with
the vanishing discs $D_0$ and $D_1$ near its extremities, it is in
general not Lagrangian (except if $S_0=S_1$, in which case we can
choose $\Psi_s=\Id$ and $S_\gamma=D_0\cup D_1$ is Lagrangian).

\begin{remark} \label{rem:iso_match}
Any path isotopic to a matching path through an isotopy fixing the
end points in $\CP^1$ and avoiding the critical points is also a
matching path. This is because the symplectic connection
associated to the pencil is Hamiltonian (i.e., parallel transport
along contractible loops generates Hamiltonian isotopies).
\end{remark}

\begin{definition}\label{def:isomatching}
Suppose that we have a family of Lefschetz pencils $\{ \phi_t
\}_{t\in [0,1]}$ and a family of paths $\gamma_t\subset\CP^1$ such that
the endpoints of $\gamma_t$ are critical values of $\phi_t$.
We say that the paths $\gamma_0$ and $\gamma_1$ are
homotopic if the arcs $\gamma_t$ pass through critical values of $\phi_t$
only for a finite number of values of the parameter $t$, and if
whenever $\gamma_{t_j}$ passes through a point
$z_j=\gamma_{t_j}(s_0)\in\mathrm{crit}(\phi_{t_j})$ with $0<s_0<1$, the
vanishing cycles $S'$ and $S''$ associated to the arcs
$\gamma_{t_j}|_{[0,s_0/2]}$ and $\gamma_{t_j}|_{[s_0/2,s_0]}$ can be made
mutually disjoint by compactly supported Hamiltonian isotopies inside
$\phi_{t_j}^{-1}(\gamma_{t_j}(s_0/2))$.
\end{definition}

Consider two homotopic paths $\gamma_0$ and $\gamma_1$, joined by
a family $\gamma_t$ as in the definition. Assume that $\gamma_0$
is a matching path: then $\gamma_1$ is also a matching path. This
can be seen by considering the families of vanishing cycles
$S_{t,0},\,S_{t,1}\subset
\Sigma_t=\phi_t^{-1}(\gamma_t(\frac12))$, $0\le t\le 1$. More
precisely, to handle the situation where $\gamma_t$ passes through
a critical value of $\phi_t$, we need to consider not just
parallel transport but also Hamiltonian isotopies $\rho_t$ inside
the fibers above $\gamma_t(\epsilon)$ and $\gamma_t(1-\epsilon)$
for some small $\epsilon>0$; $S'_{t,0}$ is then defined by taking
parallel transport along $\gamma_t|_{[0,\epsilon]}$, then applying
$\rho_t$, and then parallel transport along
$\gamma_t|_{[\epsilon,1/2]}$ (and similarly for $S'_{t,1}$).
Definition \ref{def:isomatching} implies that, by choosing
$\rho_t$ suitably, we can ensure that the perturbed vanishing
disks avoid the critical values of $\phi_t$ for all $t$, and hence
define Lagrangian spheres $S'_{t,0}$, $S'_{t,1}$ that depend
continuously on $t$ and are Hamiltonian isotopic to the vanishing
cycles for all $t$ (we can also assume $\rho_0=\rho_1=\Id$). With
this understood, by assumption $S'_{0,0}$ is isotopic to
$S'_{0,1}$, so $S'_{t,0}$ and $S'_{t,1}$ are mutually isotopic for
all $t$, and even Hamiltonian isotopic because of exactness. In
fact, we can find a family of exact symplectomorphisms identifying
the symplectic submanifolds $\Sigma_t$ with each other, and the
Hamiltonian isotopy between $S'_{1,0}$ and $S'_{1,1}$ can then be
realized by juxtaposition of the families of exact Lagrangian
spheres $(S'_{1-t,0})_{0\le t\le 1}$, $(S'_{0,s})_{0\le s\le 1}$,
and $(S'_{t,1})_{0\le t\le 1}$.

So being a matching path is a property that
depends only on the relative homotopy type of the path and not on
a particular realization.

The following Lemma is an unpublished result of Donaldson. We
thank Paul Seidel for communicating it to us. We give the
parametric version for completeness.

\begin{lemma} \label{lemm:match}
Let $\gamma_t$ be a $1$-parametric family of matching paths in a
family of symplectic pencils $\phi_t$, $t\in [0,1]$ (of large enough
degree). Then there exists a continuous family of Lagrangian spheres $S_t$
in $M$ such that each of them is smoothly isotopic to
$S_{\gamma_t}$.

Moreover, if for $j=0,1$ the vanishing cycles $S_{j,0}$ and
$S_{j,1}$ coincide over the reference fiber, then the family of
Lagrangian spheres $S_t$ can be chosen in such a way that
$S_j=S_{\gamma_j}$ for $j=0,1$.
\end{lemma}

\begin{proof}
We first assume that the isotopy does not cross critical points.
Remove a fiber $\phi_t^{-1}(z_t)$ from each of the pencils
$\phi_t$, where the point $z_t$ lies outside of the image of
$\gamma_t$. After composition with a projective automorphism of
$\CP^1=\C\cup\{\infty\}$, we can assume that $z_t=\infty$ and
$\gamma_t(1/2)=0$. The fiber above infinity is the zero set of an
asymptotically holomorphic section $s_t$, of which we may assume
without loss of generality that it is transverse to zero provided
$k$ is large enough (see \cite{AMP02} for details of how to get
transversality in a generic fiber of a Donaldson pencil). {}From
now on we will consider the restrictions of the pencils $\phi_t$
to the open manifolds $M_t$, still denoted by $\phi_t:M_t \to \C$.

We are going to perturb the symplectic form in the fibers of $\phi_t$ over
$\gamma_t|_{[1/2-\e, 1/2]}$, where $\epsilon>0$ is a small constant,
in order to change the symplectic connection and make the two vanishing
cycles match. After a reparametrization of $\gamma_t$, we can assume that
the interval over which we perturb the symplectic form is $[0,1]$ instead
of $[\frac12-\e,\frac12]$.
Identify small closed neighborhoods of $\gamma_t([0,1])$ in $\CP^1$ with
$D=[0,1]\times [-1,1]$, in such a way that $\gamma_t([0,1])$ is mapped to
$[0,1]\times \{0 \}$. Then
we can construct a family of charts
  $$
  \Phi_t: D\times F \stackrel{\sim}{\longrightarrow}\phi_t^{-1}(D) \cap M_t,
  $$
where $F$ is an open symplectic (in fact Stein) fiber of $\phi_t$,
by identifying symplectically all the fibers of $\phi_t$
above $D$ in such a way that parallel transport along $[0,1]\times \{0 \}$ is
horizontal. Hence, calling $(x,y)$ the coordinates on $D=[0,1]\times[-1,1]$,
  $$
  \Phi_t^*(\omega)= \sigma_t + \alpha_t \wedge dx + \beta_t
  \wedge dy + f_t\, dx \wedge dy,
  $$
where $\sigma_t=\sigma$ is the pullback of a symplectic form on
$F$, constant over $D$ and independent of $t$ (it follows from
Moser's argument that the fibers of the various pencils $\phi_t$
are mutually symplectomorphic), while $\alpha_t(x,y)$ and
$\beta_t(x,y)$ are $1$-forms in the fibers, with $\a_t(x,0)\equiv
0$, and $f_t(x,y)$ is a positive function. The closedness of
$\omega$ imposes the following relations:
  $$
  d\alpha_t=0, \quad d\beta_t=0, \quad df_t =
  \frac{\bd \beta_t}{\bd x}-\frac{\bd \alpha_t}{\bd y},
  $$
where the exterior differentials are only applied to the fiber
directions (the $x$ and $y$ coordinates are treated as
parameters). Choose a vector $v=(v_1,v_2)\in TD$. The unique horizontal
lifting of $v$ is a vector $(v, X_v)$, where $X_v\in TF$ is determined by
the equation
  $$
  i_{X_v} \sigma= v_1\,\alpha_t + v_2\,\beta_t.
  $$
In particular the lifting of the segment $[0,1]\times \{ 0 \}$ gives
rise to a vector field with $X_v=0$.

Let $\Psi_{t,s}$ be a compactly supported Hamiltonian isotopy of the fiber $F$
sending the vanishing sphere $S_{t,0}$ to $S_{t,1}$ (recall that the
variable $t$ parametrizes the family and $s$ is the time
parameter for the Hamiltonian flow of each member of the family).
We want to change the symplectic connection (changing the
symplectic form) in such a way that the horizontal lift of the segment
$[0,1]\times\{0\}$ generates $\Psi_{t,s}$ instead of the identity map.
Let $H_{t,s}$ be a family of Hamiltonian functions generating
$\Psi_{t,s}$. Without loss of generality we may assume that $H_{t,s}$
vanishes identically for all $t\in [0, \d]\cup [1-\d,1]$, for some small
$\d>0$. We define a function over $[\frac14,\frac34]\times \{ 0 \} \times F
\subset D\times F$
by the formula $\hat{F}_t(x,0, p)=H_{t,2(x-1/4)}(p)$. Extend $\hat{F}_t$ to
all of  $D\times F$ in such a way that it vanishes outside of $[\frac14,
\frac34]\times [-\frac12,\frac12]\times F$. Define $\alpha_t'=\alpha_t+ d\hat{F}_t$
and $f_t'= f_t- \frac{\bd \hat{F}_t}{\bd y}$, and consider the new 2-form
  $$
  \omega_t'=\sigma + \alpha_t' \wedge dx + \beta_t \wedge dy +
  f_t'\, dx \wedge dy.
  $$
The closed 2-form $\omega_t'=\Phi_t^*(\omega)+d(\hat{F}_t\,dx)$
coincides with $\Phi_t^*(\omega)$ in a neighborhood of the
boundary of $F\times D$, therefore we can extend $\omega_t'$ to a
closed 2-form over $M$, coinciding with $\omega$ outside of
$F\times D$. Moreover, by construction the horizontal lift of the
segment $[0,1]\times \{0\}$ with respect to the symplectic
connection induced by $\omega'_t$ generates the flow $\Psi_{t,s}$,
so that the vanishing cycles now coincide. In general, although
$\omega'_t$ is vertically non-degenerate, it need not be
symplectic. We construct a new family of forms
 $$
 \tilde{\omega}_t= \omega_t' + \phi_t^*(C_t \, dx\wedge dy),
 $$
where $C_t$ is a positive real constant. For $C_t$ large enough
the family $\tilde{\omega}_t$ is symplectic. Moreover, if
$\Psi_{t,0}=\Psi_{t,1}=\Id$, then $\omega'_0 = \omega'_1 = \omega$
and the continuous family of constants $C_t$ can be chosen to
satisfy $C_0=C_1=0$. Recall that $\tilde{\omega}_t$ and
$\omega_t'$ generate the same symplectic connection and therefore
the matching paths $\gamma_t$ generate Lagrangian spheres
$\tilde{S}_t$ for the symplectic structure $\tilde{\omega}_t$.

There is a biparametric family of exact symplectic structures
$\tilde{\omega}_{t,s}$ such that $\tilde{\omega}_{t,0}=
\tilde{\omega}_t$ and $\tilde{\omega}_{t,1}=\omega$. To construct
it we first shrink to zero the perturbation in $\alpha_t'$ and
then we shrink to zero the constant $C_t$ (the order is important
to keep all the forms non-degenerate). Hence, each of these forms
can be written as
$\tilde{\omega}_{t,s}=\omega+a_{t,s}\,d(\hat{F}_t\,dx)+
C_{t,s}\,\phi_t^*(dx\wedge dy)$, for some $a_{t,s},\,C_{t,s}\ge
0$. By applying Moser's trick to the family of forms $\{
\tilde{\omega}_{t,s} \}_{s\in[0,1]}$, we can find a family of
diffeomorphisms $\psi_{t,s',s}$, defined over open subsets of
$M_t$, such that $\psi_{t,s', s}(\tilde{\omega}_{t,s}) =
\tilde{\omega}_{t,s+s'}$. In the case where
$\Psi_{0,s}=\Psi_{1,s}=\Id$, we have
$\psi_{0,s',s}=\psi_{1,s',s}=\Id$. Now, we want to push-forward
$\tilde{S}_t$ using the flow $\psi_{t,s',s}$ in order to obtain a
family of Lagrangian spheres for the initial symplectic structure
$\omega$. However, since $M_t$ is an open manifold,
$\psi_{t,s',s}$ is not necessarily well defined everywhere.

This difficulty can be avoided by using the exactness of the
symplectic forms $\tilde{\omega}_{t,s}$ on $M_t$, and the
existence of a Liouville vector field $v_{t,s}$ transverse to the
boundary. More precisely, recall that over $M_t$ we can write
$\omega=d\theta_t$, where $-k\theta_t$ is the imaginary part of
the connection 1-form of $L^{\otimes k}$ in the trivialization
given by the defining section $s_t$ of the fiber at infinity. In
other terms, $\theta_t=-\frac{1}{k}\mathrm{Im}(s_t^{-1}\nabla
s_t)$. The vector field $v_{t}$ such that
$\omega(v_t,\cdot)=\theta_t$ is then a Liouville vector field for
$\omega$, and since $s_t$ is asymptotically holomorphic and
transverse to $0$, near the boundary of $M_t$ this vector field
points everywhere outwards -- in fact, near the boundary $v_t$
coincides up to $O(k^{-1/2})$ with the gradient vector field of
$-\frac{1}{k}\log |s_t|$. When we perturb the symplectic form to
$\tilde{\omega}_{t,s}$ the existence of a Liouville vector field
is preserved: we can write
$\tilde{\omega}_{t,s}=d\tilde{\theta}_{t,s}$, where
$$\tilde{\theta}_{t,s}=\theta_t+a_{t,s}\,\hat{F}_t\,dx+{\tfrac{1}{2}}C_{t,s}\,
\phi_t^*(x\,dy-y\,dx).$$ The corresponding Liouville vector field
$\tilde{v}_{t,s}$ (defined by
$\tilde{\omega}_{t,s}(\tilde{v}_{t,s},\cdot)=
\tilde{\theta}_{t,s}$) is still pointing outwards near the
boundary of $M_t$, as can be seen by checking that
$\tilde{\theta}_{t,s}\wedge \tilde{\omega}_{t,s}$ has the required
positivity property (here one uses the fact that the first
perturbation term $\hat{F}_t\,dx$ is supported inside a compact
subset of $M_t$ and, more importantly, the non-negativity of
$C_{t,s}$).

 Denote $k_{t,s,\tau}$ the flow generated by integrating $-\tilde{v}_{t,s}$
over a time $\tau$; this flow is well-defined everywhere and pushes
inwards near the boundary of $M_t$. It is easy to construct
two families of open sets $D_{t,1}\subset D_{t,2}\subset M_t$ with the
following properties:
\begin{enumerate}
\item $\tilde{S}_t \subset D_{t,1}$. \item There is a small $\e>0$
(independent of $t$ and $s$ by compactness) for which $\psi_{t,
\e,s}(D_{t,1}) \subset D_{t,2}$; assume $l=\frac{1}{\e}\in \N$.
\item There is a $\lambda_t>0$ such that
$k_{t,s,\lambda_t}(D_{t,2})\subset D_{t,1}$ for all $s,t$.
\end{enumerate}
In the case where $\Psi_{0,s}=\Psi_{1,s}=\Id$, we can choose
$D_{j,1}=D_{j,2}$ and $\lambda_j=0$ for $j\in\{0,1\}$.

Observe that the flow $k_{t,s,\lambda_t}$ preserves the property
of being Lagrangian with respect to $\tilde{\omega}_{t,s}$.
Therefore, we can proceed as follows. Start with $\tilde{S}_{t,0}=
\tilde{S}_{t}$, and let $\tilde{S}_{t,1}= k_{t, \epsilon,
\lambda_t} \circ \psi_{t, \e, 0} (\tilde{S}_{t,0})$. By
construction $\tilde{S}_{t,1} \subset D_{t,1}$. We can then repeat
the process, defining
  $$
  \tilde{S}_{t,j+1}= k_{t, (j+1)\epsilon, \lambda_t} \circ
  \psi_{t, \epsilon, j\epsilon} (\tilde{S}_{t,j}).
  $$
At the end we get a family of Lagrangian spheres $\tilde{S}_{t,l}$
with respect to the initial symplectic structure $\omega$, which
completes the proof in the case where the isotopy does not cross
any critical point of $\phi_t$.

In the general case, the argument remains the same.
The only important observation is that, by definition, suitably chosen Hamiltonian
isotopies in the fibers can be used to ensure that the spheres we consider
stay away from the vanishing cycles at the critical points hit
by the family of matching paths $\gamma_t$.
Hence, we can still find a family of Hamiltonian isotopies $\Psi_{t,s}$ in the fibers
of $\phi_t$, sending the vanishing sphere $S_{t,0}$ to $S_{t,1}$, and with
support in a compact subset of $M_t$ disjoint from the support of the
generalized Dehn twists arising as monodromies around the different
critical values encountered by $\gamma_t$.
\end{proof}

\begin{remark} \label{rem:unique_sph}
 The previous construction does not provide a canonical sphere
 associated to the matching path. However, all the possible choices
 are Lagrangian (and hence Hamiltonian) isotopic. This is because
 the spaces of choices that appear in the proof are all
 path connected.
\end{remark}

\section{Isotopies of Lagrangian submanifolds.}
\label{sec:isotopies}

The above results allow us to tentatively identify Lagrangian
spheres with matching paths. By Theorem \ref{thm:main} we know
that, given any Lagrangian sphere, there exists a pencil for which
it fibers above a matching path. Conversely, by Lemma
\ref{lemm:match} we know that a matching path gives rise to a
Lagrangian sphere up to isotopy. It is therefore natural to look
for a general result identifying homotopy classes of matching
paths and Lagrangian spheres. For this purpose, we now study more
carefully the behavior of our construction in the presence of
isotopic Lagrangian spheres (which will correspond to homotopic
paths, cf.\ Remark \ref{rem:iso_match}).

Consider a fixed sequence of Donaldson pencils $\phi_k^0$, and introduce the
following
\begin{definition}\label{def:asymptstable}
A sequence of matching paths $\gamma_k$ for the
sequence of Donaldson pencils $\phi_k^0$, associated to Lagrangian
spheres in the same Hamiltonian isotopy class for all $k$, is called
asymptotically stable if for large $k$ there exists a
family $(\phi_{k,t})_{0\le t\le 1}$ of Donaldson pencils with $\phi_k^0=
\phi_{k,0}$, such that the path $\gamma_k$ is homotopic through the family
$\phi_{k,t}$ to the path
\begin{eqnarray*}
\hat{\gamma}_k: [0,1] & \to & \CP^1 \\
t & \to & e^{2ik^{1/2}t+k^{-1/2}t}
\end{eqnarray*}
and the Lagrangian sphere associated to the matching path $\hat{\gamma}_k$
arises from the construction described in the proof of Theorem
\ref{thm:main}
(this implies in particular a number of transversality conditions on
$\phi_{k,1}$ near the Lagrangian sphere).
\end{definition}

\begin{remark}
The sequences generated by Theorem \ref{thm:main} are asymptotically
stable sequences once we isotop the resulting family of pencils
$\phi_k$ to the fixed family $\phi_k^0$. In fact, these are the
only computable examples.
\end{remark}

The following result is a direct corollary of Theorem \ref{thm:main_par}
and Definition \ref{def:asymptstable}:

\begin{theorem} \label{thm:main_para}
Let $S_0$ and $S_1$ be two Hamiltonian isotopic Lagrangian spheres,
associated to asymptotically stable sequences of
matching paths $\gamma_{k,0}$ and $\gamma_{k,1}$ in the pencils $\phi_k^0$.
Then, for large enough $k$, there exists a family of Lefschetz pencils
$(\phi_{k,t})_{0\le t\le 1}$, with $\phi_{k,0}=\phi_{k,1}=\phi_k^0$,
such that the matching path $\gamma_{k,0}$ is homotopic to
$\gamma_{k,1}$ via the family $\phi_{k,t}$.
\end{theorem}

Define the set
 $$
 \SM= \{ [\gamma_k]: \gamma_k ~ \mbox{
 is an asymptotically stable sequence for} ~ \phi_k^0,
 ~ k\ \mbox{large enough} \},
 $$
where $[\gamma_k]$ denotes the homotopy class of the matching path
$\gamma_k$. Define the set
 $$
 Pen = \{ \{\phi_k\}: \phi_k\ \mbox{is ~ a ~ sequence ~ of  ~ Donaldson ~
 pencils} \}.
 $$
There is a natural action
 \begin{equation}
 \pi_1(Pen, \{\phi_k^0\}) \times \SM \to \SM,
 \label{eq:loop_pen}
 \end{equation}
defined by transport along a family of pencils. Then we have:
\begin{theorem} \label{thm:equiv}
The set of Hamiltonian isotopy classes of Lagrangian spheres in
the symplectic manifold $(M,\omega)$ is in one-to-one
correspondence with the set $\SM / \pi_1(Pen,\{\phi_k^0\})$.
\end{theorem}

In principle this result reduces (at least in dimension $4$) the
problem of the classification of Lagrangian spheres in a
symplectic manifold to the purely algebro-combinatorial problem of
classifying matching paths in a sequence of pencils. However,
the notion of asymptotic stability of matching paths is a very unnatural
one, and makes things much less practical.

To eliminate this requirement, and to simplify the discussion, one
should consider the behavior of pencils and matching paths under
degree doubling, i.e.\ upon passing from sections of $L^{\otimes
k}$ to sections of $L^{\otimes 2k}$. In general this requires an
understanding of the ``degree doubling'' process (see e.g.\
\cite{AKdoubling} \cite{Smdoubling}). For simplicity we restrict
ourselves to the algebraic case, where estimated transversality is
not needed. Recall that, in that case, there is a natural way of
inducing a sequence of matching paths in a sequence $\phi_k^0$ of
pencils ($k=2^l k_0,\ l\in \N$), built as follows: starting from a
given pencil $\phi^0_{k}=s_{k}^1/s_{k}^2$, construct a new pencil
$\phi^0_{2k}=s_{2k}^1/s_{2k}^2$ whose defining sections are
arbitrarily small generic perturbations of
$\sigma_{2k}^1=s_{k}^1\otimes s_{k}^1$ and
$\sigma_{2k}^2=s_{k}^1\otimes s_{k}^2$. Repeating the process we
get a sequence of pencils $\phi_{k}^0$, $k=k_0 2^l$. It is easy to
check that the set of critical values of $\phi_k^0$ identifies
naturally with a subset of the critical values of $\phi_{2k}^0$.
More precisely, there is a diffeomorphism from an open ball
$B(0,R)$ of $\CP^1$ to itself which takes the critical values of
$\phi_k^0$ to critical values of $\phi_{2k}^0$, and such that the
rest of the critical values of $\phi_{2k}^0$ remain outside that
ball; in fact, the monodromy of $\phi_k^0$ naturally ``embeds''
into that of $\phi_{2k}^0$ \cite{AKdoubling} \cite{Smdoubling}.
This makes it possible to build from a matching path
$\gamma_{k_0}$ for $\phi_{k_0}^0$ a sequence of matching paths
$\gamma_k$ for $\phi_k^0$, for all $k=2^l k_0$. It is easy to
check that the Lagrangian sphere associated to each element of the
family is always the same up to Hamiltonian isotopy. We will call
``natural sequence'' the sequence of matching paths obtained in
this manner from a given matching path. Then Theorem
\ref{thm:equiv} can be reformulated as

\begin{theorem} \label{thm:algeb}
The set of natural sequences $\phi_k^0$ in a family of Lefschetz
pencils on a projective manifold up to the action of $\pi_1(Pen,
\{\phi_k^0\})$ is in bijection with the set of Hamiltonian isotopy
classes of Lagrangian spheres.
\end{theorem}

\section{Pencil automorphisms, symplectomorphisms, and matching paths}
\label{sec:automorphisms}

In this section we discuss the connection between automorphisms
of the monodromy data of a Lefschetz pencil and symplectomorphisms of the
total space of the pencil; most of the discussion follows ideas of
Donaldson and Seidel (see in particular Section (1d) of \cite{Sei03}).

\subsection{The group of pencil automorphisms}
To any Lefschetz pencil one can attach a group of ``pencil automorphisms'',
which can be viewed either geometrically (a pencil automorphism is then
a homeomorphism which lifts to a diffeomorphism of the corresponding
Lefschetz fibration on the blown-up manifold, mapping fibers to fibers and inducing fiberwise
symplectomorphisms), or (taking isotopy classes)
combinatorially in terms of monodromy data. We will adopt the
combinatorial point of view here.

Recall that the monodromy of a symplectic
Lefschetz pencil $\phi:M-N\to \CP^1$ is defined by fixing a base point
$\alpha_0\in\CP^1$ and considering the isotopy classes of the
symplectomorphisms of the reference fiber
$\Sigma_0=\overline{\phi^{-1}(\alpha_0)}$ induced by
parallel transport along loops in $\CP^1-\mathrm{crit}(\phi)$.
The monodromy morphism takes values in the symplectic
mapping class group of the fiber relatively to the base points. More precisely, let
$$\mathrm{Map}(\Sigma_0,N)=\pi_0\{g\in\mathrm{Symp}(\Sigma_0,
\omega_{|\Sigma_0})\ |\ \forall p\in N,\ g(p)=p\ \mathrm{and} \
dg(p)=\mathrm{Id}\}.$$ Then, after removing a neighborhood of a
smooth fiber $\Sigma_\infty$ of the pencil, we obtain a fibration
over a large disc $\Delta$ containing all the critical values of
$\phi$, whose monodromy is given by a morphism
$$\theta:\pi_1(\Delta-\mathrm{crit}(\phi),\alpha_0)\to \mathrm{Map}(\Sigma_0,N).$$

Let $r$ be the number of critical points of $\phi$, which we assume
to all lie in distinct fibers, and let $\mathcal{A}=\mathrm{crit}(\phi)
\subset\Delta\subset\CP^1$ be the
set of critical values. For simplicity we assume that the base point
$\alpha_0$ lies on the boundary of $\Delta$. Choosing a set of geometric
generators $\gamma_1,\dots,\gamma_r$ of
$\pi_1(\Delta-\mathcal{A},\alpha_0)$ (each encircling one of the
critical values), each
$\theta(\gamma_i)$ is a Dehn twist about an (exact) Lagrangian sphere
$L_i\subset \Sigma_0-N$ (the vanishing cycle associated to $\gamma_i$).
Any loop in $\pi_1(\Delta-\mathcal{A},\alpha_0)$ that bounds a disc
containing exactly one critical value can be written in the form
$g^{-1}\gamma_i g$, and is mapped by $\theta$ to a Dehn twist about the
Lagrangian sphere $\theta(g)(L_i)$.

\begin{definition}
Let $\hat{\mathcal{M}}$ be the set of all pairs $(\zeta,L_\zeta)$ where
$\zeta:\pi_1(\Delta-\mathcal{A},\alpha_0)\to\mathrm{Map}(\Sigma_0,N)$ is a group
homomorphism and $L_\zeta$ is a map from the set $\Pi$ of all conjugates of
geometric generators in $\pi_1(\Delta-\mathcal{A},\alpha_0)$ to the set of
Hamiltonian isotopy classes of (exact) Lagrangian spheres in $\Sigma_0-N$,
such that: $(i)$ $\forall \gamma\in\Pi$, $\zeta(\gamma)$ is the Dehn twist
about $L_\zeta(\gamma)$; $(ii)$ $\forall\gamma\in\Pi$,
$\forall g\in \pi_1(\Delta-\mathcal{A},\alpha_0)$,
$L_\zeta(g^{-1}\gamma g)=\zeta(g)(L_\zeta(\gamma))$.
\end{definition}

Then we can associate to the Lefschetz pencil $\phi$ an element
$\hat\theta=(\theta,L_\theta)\in\hat{\mathcal{M}}$, where $\theta$ is
the monodromy morphism and $L_\theta$ is characterized by the property that
$L_\theta(\gamma_i)=[L_i]$. When $\dim M=4$, one can recover the vanishing
cycles from the monodromy morphism (using the exactness property to
determine the Hamiltonian isotopy class), so in that case $\hat\theta$
contains no more information than $\theta$; in higher dimensions it is
unknown whether $\theta$ determines $\hat\theta$.

It is a result of Gompf that the ``enhanced'' monodromy morphism
$\hat\theta$ determines the topology of the Lefschetz pencil over
a large disc containing all the critical values of $\phi$, and up
to a choice in $\pi_1\mathrm{Symp}(\Sigma_0,N)$ (describing the
attaching map near the fiber at infinity), the symplectic manifold
$(M,\omega)$ up to symplectic isotopy \cite{Gompf}. This result
relies on the fact that the fibers of the pencil are Poincar\'e
dual to a multiple of $[\omega]$: in general, applying Thurston's
argument to a Lefschetz fibration yields symplectic forms that are
only canonical up to deformation equivalence, but in the case of a
pencil, after blowing down the exceptional divisor all these forms
become isotopic up to scaling by a constant factor.

We will also consider the group
 $$
 \mathcal{B}(\Sigma_0,N)=\pi_0\{g\in\mathrm{Symp}(\Sigma_0,\omega_{|\Sigma_0})\ |
 \ g(N)=N\}.
 $$
This group acts on $\mathrm{Map}(\Sigma_0,N)$ by conjugation, and given an
element $g\in\mathcal{B}(\Sigma_0,N)$ we denote by $g_*$ the corresponding
automorphism of $\mathrm{Map}(\Sigma_0,N)$. Combining this with the
natural action of $\mathcal{B}(\Sigma_0,N)$ on the set of isotopy classes
of Lagrangian spheres, we obtain an action of $\mathcal{B}(\Sigma_0,N)$
on $\hat{\mathcal{M}}$, given by $g_*: (\zeta,L_\zeta)\mapsto(g_*\circ \zeta,g\circ
L_\zeta)$.

When $M$ is a 4-manifold, let $h$ be the genus of $\Sigma_0$ and $n$
the number of base points: then
$\mathrm{Map}(\Sigma_0,N)=\mathrm{Map}_{h,n}$ is the mapping class
group of a genus $h$ surface with $n$ boundary components, and
$\mathcal{B}(\Sigma_0,N)=\mathcal{B}_{h,n}$ is the braid group on
$n$ strings on a genus $h$ surface.

Finally, denote by $B_r$ the classical braid group on $r$
strings, which can be viewed either as the fundamental group of the
configuration space of $r$ distinct points in the disc
($B_r=\pi_1(\mathrm{Conf}_r(\Delta),\mathcal{A})$) or in terms of isotopy classes
of compactly supported orientation-preserving diffeomorphisms of $\Delta$
mapping $\mathcal{A}$ to itself. The latter description lets us associate
to any braid $b\in B_r$ an automorphism $b_*$ of $\pi_1(\Delta-\mathcal{A},\alpha_0)$.
Hence there is a natural right action of $B_r$ on $\hat{\mathcal{M}}$ (by composition).
With the above notations, we can make the following definition:

\begin{definition}
The group of combinatorial automorphisms of $\phi$ is $$\Gamma(\phi)=
\{(b,g)\in B_r\times \mathcal{B}(\Sigma_0,N)\ |\ \hat\theta\circ b_*=
g_*(\hat\theta)\}.$$
\end{definition}

In other words, $\Gamma(\phi)$ is the stabilizer of the monodromy data
$\hat\theta\in\hat{\mathcal{M}}$ for the natural actions of the braid group (by {\it Hurwitz
moves}) and automorphisms of the fiber (by conjugation).

There is a natural homomorphism $\rho:\Gamma(\phi)\to
\pi_0\mathrm{Symp}(M,\omega)$, which can be described as follows.
Given an element $(b,g)\in\Gamma(\phi)$, choose a geometric
representative of the braid $b$ inside a cylinder $\Delta\times
[0,1] \subset\CP^1\times [0,1]$, and use it to build a
one-parameter family of Lefschetz pencils $\phi_t$ (each having
the same total space $M$ and the same monodromy, but with critical
values depending on $t\in [0,1]$ as prescribed by the braid $b$).
By identifying the two ends of the cylinder we can close the braid
$b$ and obtain a link in $\Delta\times S^1\subset\CP^1\times S^1$.
The equality $\hat\theta\circ b_*=g_*(\hat\theta)$ makes it
possible to identify the Lefschetz pencils $\phi_0=\phi$ and
$\phi_1$ by means of a fiberwise symplectomorphism in the isotopy
class $g\in\mathcal{B}(\Sigma_0,N)$, up to a symplectic isotopy
(using Gompf's characterization of symplectic structures on
Lefschetz pencils~\cite{Gompf}). This yields a family of
symplectic Lefschetz pencils parametrized by elements of $S^1$,
whose total space $W$ carries a structure of symplectic fiber
bundle $(M,\omega)\to W\to S^1$; its monodromy is the element of
$\pi_0\mathrm{Symp}(M,\omega)$ naturally associated to the
automorphism~$(b,g)$.

The following asymptotic surjectivity result is a direct consequence of
Donaldson's result of uniqueness up to isotopy \cite{Do00}:

\begin{proposition}\label{prop:rhosurj}
Let $\phi_k:M-N_k\to\CP^1$, $k\gg 0$ be a sequence of Donaldson's
symplectic Lefschetz pencils. Then for every $\eta\in
\pi_0\mathrm{Symp}(M,\omega)$ there exists an integer $k(\eta)$
such that $\eta$ belongs to the image of the natural homomorphism
$\rho: \Gamma(\phi_k)\to\pi_0\mathrm{Symp}(M,\omega)$ for all
$k\ge k(\eta)$.
\end{proposition}

\begin{proof}
Let $\tilde\eta$ be a symplectomorphism in the isotopy class
$\eta$, and equip $M$ with $\omega$-compatible almost-complex
structures $(J_t)_{t\in [0,1]}$ such that $J_0=\tilde\eta_*(J_1)$.
By identifying the boundaries of $M\times [0,1]$ via the
symplectomorphism $\tilde\eta$, we can build a symplectic fiber
bundle over $S^1$ with fiber $(M,\omega)$ and monodromy $\eta$.
Donaldson's construction of Lefschetz pencils \cite{Do00} applies
to this one-parameter family and provides, for $k\gg 0$, a family
of symplectic Lefschetz pencils
$\tilde\phi_{k,t}:M-N_{k,t}\to\CP^1$ such that
$\tilde\phi_{k,1}=\tilde\phi_{k,0}\circ \tilde\eta$ (in
\cite{Do00} the case where the parameter is in $t\in [0,1]$ is
used to prove uniqueness up to isotopy; the case where the
parameter belongs to $S^1$ is exactly identical). We can also
easily require the critical values of $\tilde\phi_{k,t}$ to remain
distinct for all values of $t$. The monodromy of the family
$\tilde\phi_{k,t}$ can then naturally be expressed as an element
$(b,g)\in \Gamma(\tilde\phi_{k,0})$ by considering the motion of
the critical values as $t$ varies in $S^1$ (which gives the braid
$b$) and the isotopy class of the induced symplectomorphism of a
generic fiber (which gives $g$); by construction we have
$\rho(b,g)=\eta$. Finally, it is known that the Donaldson pencils
$\tilde\phi_{k,0}$ and $\phi_k$ are mutually isotopic for large
enough $k$ \cite{Do00}, so that we can naturally identify
$\Gamma(\tilde\phi_{k,0})$ with $\Gamma(\phi_k)$.
\end{proof}

It is an interesting question to ask whether there is a value of $k$ for
which the morphism $\rho$ is surjective. A positive answer might follow from
a better understanding of the behavior of the group of pencil automorphisms
under stabilization by degree doubling.

\subsection{Matching paths, Dehn twists, and trivial automorphisms}

Matching paths provide an explicit way to view Dehn twists
along Lagrangian spheres in $M$ in this context. We begin with
some terminology. Let $\delta:[0,1]\to\Delta\subset\CP^1$ be an embedded
arc meeting $\mathcal{A}$ only at its endpoints $\delta(0)$ and
$\delta(1)$. Then the positive {\it half-twist} along $\delta$ is
the braid $\sigma_\delta\in B_r$ which exchanges the two points
$\delta(0)$ and $\delta(1)$ by a 180-degree counterclockwise rotation
in a small neighborhood of $\delta([0,1])$. Also, choose an arc $\eta$
in $\Delta-(\mathcal{A}\cup\delta)$ joining the base point $\alpha_0$ to
the point $\delta(0)$; let $\eta'$ be the oriented boundary
of a small tubular neighborhood of $\eta$, and let $\eta''$ be the image
of $\eta'$ by the half-twist $\sigma_\delta$. Then $\eta',\eta''\in\Pi$
are conjugates of geometric generators of $\pi_1(\Delta-\mathcal{A},\alpha_0)$;
the pair $(\eta',\eta'')$ is called a {\it supporting pair} for the arc
$\delta$, and is unique up to simultaneous conjugation
$(\eta',\eta'')\mapsto (g^{-1}\eta'g,g^{-1}\eta''g)$
in $\pi_1(\Delta-\mathcal{A},\alpha_0)$. In particular, the pair of vanishing
cycles $(L_\theta(\eta'),L_\theta(\eta''))$ is uniquely determined up to
simultaneous action of an automorphism in the monodromy group
$(L_\theta(\eta'),L_\theta(\eta''))\mapsto (\theta(g)(L_\theta(\eta')),
\theta(g)(L_\theta(\eta'')))$.

\begin{proposition}
Let $\gamma:[0,1]\to\Delta\subset\CP^1$ be a matching path for the Lefschetz
pencil $\phi$, corresponding to a Lagrangian sphere $S_\gamma\subset M$.
Then $(\sigma_\gamma,1)\in \Gamma(\phi)$, and
$\rho(\sigma_\gamma,1)=[\tau_{S_\gamma}]$ $($where $\tau_{S_\gamma}$ is
the Dehn twist along $S_\gamma)$.
\end{proposition}

\begin{proof}
Let $(\eta',\eta'')$ be a supporting pair for $\gamma$. By choosing
geometric generators for the complement of $(\eta\cup\gamma)$ in $\Delta$, we can
complete the pair $(\eta',\eta'')$ to an ordered collection
of generators $(\gamma_i)_{1\le i\le r}$ of $\pi_1(\Delta-\mathcal{A},\alpha_0)$
such that $\gamma_i\cap\gamma=\emptyset$ for all $i\ge 3$.
Then the action of $\sigma_\gamma$ on
$\pi_1(\Delta-\mathcal{A},\alpha_0)$ is given by $\sigma_\gamma(\eta')=\eta''$,
$\sigma_\gamma(\eta'')=\eta''\eta'(\eta'')^{-1}$, and
$\sigma_\gamma(\gamma_i)=\gamma_i$ for all $i\ge 3$.

By definition of
a matching path, the vanishing cycles associated to $\eta'$ and $\eta''$
are Hamiltonian isotopic, i.e.\ $L_\theta(\eta')=L_\theta(\eta'')$. This
implies immediately that
$\hat\theta\circ (\sigma_\gamma)_*=\hat\theta$, i.e.\ $(\sigma_\gamma,1)$
is an automorphism of the pencil.

Construct a family of pencils $(\phi_t)$ whose monodromy is
$g=\rho(\sigma_\gamma,1)$ by choosing a representative of the
braid $\sigma_\gamma$ in $\Delta\times [0,1]$ and identifying
the two ends of the cylinder $t=0$ and $t=1$, as explained above.
Because the braid $\sigma_\gamma$ has a representative supported in
a small neighborhood $U_\gamma$ of the matching path $\gamma$, and because
$\phi$ is trivial over $U_\gamma$ away from the
vanishing cycle, it is easy to see that the isotopy class $g$ admits a
representative which coincides with $\mathrm{Id}$ outside of a small
neighborhood of the sphere $S_\gamma$. The argument can therefore be completed by
considering a universal local model for the family of pencils $\phi_t$
over a neighborhood of $S_\gamma$.

Consider the map $F=(F_1,F_2):\C^{n+1}\to\C^2$ defined by $$F(z_1,\dots,z_{n+1})=
(z_1^2+\dots+z_{n+1}^2,z_{n+1}).$$ For $t=\epsilon e^{i\theta}$, the
restriction of $F$ to the hypersurface $X_t=F_1^{-1}(t)=\{z_1^2+\dots+z_{n+1}^2=t\}$
induces a Lefschetz fibration $F_2:X_t\to \C$, whose generic fiber
$F_2^{-1}(u)$ is a smooth quadric in $\C^n\times\{u\}$, defined
by the equation $z_1^2+\dots+z_n^2=t-u^2$. There
are two singular fibers, and the corresponding critical values
are the two square roots of $t$. The straight line segment between the
two critical values $\pm t^{1/2}$ is a matching path
for $F_2$, and the corresponding Lagrangian sphere is
$S_t=X_t\cap (e^{i\theta/2}\R)^{n+1}$. As $\theta$ varies from $0$ to
$2\pi$, the two critical values of $F_2$ are exchanged by a half-twist
$\sigma$ along the matching path.

On the other hand, $F_1:\C^{n+1}\to\C$ provides
a local model for a neighborhood of a critical point of a Lefschetz fibration
in complex dimension $n+1$. Therefore, the monodromy of the family of symplectic
manifolds $X_{t=\epsilon\exp(i\theta)}$ as $\theta$ varies from $0$ to
$2\pi$ is the Dehn twist along the vanishing cycle associated to the
singular fiber $X_0$ of $F_1$; this vanishing cycle is precisely the
Lagrangian sphere $X_\epsilon\cap \R^{n+1}=S_\epsilon$, and so we
have $\rho(\sigma,1)=[\tau_{S_\epsilon}]$.
\end{proof}

It is also possible to construct explicitly many elements in the
kernel of the morphism
$\rho:\Gamma(\phi)\to\pi_0\mathrm{Symp}(M,\omega)$ (i.e.,
``trivial'' pencil automorphisms). Let $\gamma:[0,1]\to\Delta$ be
an embedded arc with endpoints in $\mathcal{A}$, choose a
supporting pair $(\eta',\eta'')$, and let $S',S''\subset\Sigma_0$
be Lagrangian spheres in the isotopy classes $L_\theta(\eta')$ and
$L_\theta(\eta'')$.

\begin{proposition}\label{prop:kerrho}
$(a)$~If $S'$ and $S''$ are disjoint, then
$(\sigma_\gamma^2,1)\in \Gamma(\phi)$ and \hbox{$\rho(\sigma_\gamma^2,1)=1$.}
\\
$(b)$~If $S'$ and $S''$ intersect transversely in a single
point, then $(\sigma_\gamma^3,1)\in\Gamma(\phi)$ and
$\rho(\sigma_\gamma^3,1)=1$.
\end{proposition}

\begin{proof}
The structure of the argument is the same as for the previous proposition.
We start by completing $(\eta',\eta'')$ to an ordered collection of
generators $(\gamma_i)_{1\le i\le r}$ of
$\pi_1(\Delta-\mathcal{A},\alpha_0)$ such that $\gamma_i\cap
\gamma=\emptyset$ for $i\ge 3$.
In the case where $S'$ and $S''$ are disjoint, observe that
$\sigma_\gamma^2$ maps $\eta'$ to $\tilde\eta'=\eta''\eta'(\eta'')^{-1}$, and
$\eta''$ to $\tilde\eta''=\eta''\eta'\eta''(\eta')^{-1}(\eta'')^{-1}$,
while the other generators are preserved. Since $\theta(\eta'')=[\tau_{S''}]$
has a representative supported in a neighborhood of $S''$, it acts trivially
on $S'$. Therefore $L_\theta(\tilde\eta')=\theta(\eta'')^{-1}(L_\theta(\eta'))
=L_\theta(\eta')$, and similarly $L_\theta(\tilde\eta'')=
\theta(\eta'')^{-1}(\theta(\eta')^{-1}(L_\theta(\eta'')))=L_\theta(\eta'')$,
which implies that $(\sigma_\gamma^2,1)\in\Gamma(\phi)$.

To see that the isotopy class $\rho(\sigma_\gamma^2,1)$ is trivial, observe
that the automorphism $(\sigma_\gamma^2,1)$ corresponds to the monodromy of
a $S^1$-family of pencils where two of the critical values
simply move around each other. Because the two vanishing cycles $S'$ and
$S''$ are mutually disjoint, this $S^1$-family of pencils bounds a
$D^2$-family of symplectic Lefschetz pencils on $M$ (one of which has two
critical points in the same fiber, but this is not a problem since the
vanishing cycles are disjoint). This implies that the monodromy is
trivial, i.e.\ $\rho(\sigma_\gamma^2,1)=1$.

In the case where $S'$ and $S''$ intersect transversely in a single point,
we consider the action of $\sigma_\gamma^3$ on $\pi_1(\Delta-\mathcal{A},
\alpha_0)$, which maps $\eta'$ to $(\eta''\eta')\eta''(\eta''\eta')^{-1}$
and $\eta''$ to $(\eta''\eta'\eta'')\eta'(\eta''\eta'\eta'')^{-1}$. The fact
that $(\sigma_\gamma^3,1)$ belongs to $\Gamma(\phi)$ then follows directly
from the observation that $\tau_{S'}\tau_{S''}(S')$ is Hamiltonian isotopic to
$S''$ and $\tau_{S''}\tau_{S'}(S'')$ is Hamiltonian isotopic to $S'$
(the geometric property underlying the braid relation
$\tau_{S'}\tau_{S''}\tau_{S'}\sim\tau_{S''}\tau_{S'}\tau_{S''}$).

The isotopy class $\rho(\sigma_\gamma^3,1)$ admits a representative with
support contained in the preimage of a small neighborhood $U_\gamma$ of
$\gamma$, and because $\phi$ is trivial over $U_\gamma$ away from the
vanishing cycles, we can again consider a universal local model for a
neighborhood of the configuration of vanishing cycles.
Consider the map $F=(F_1,F_2):\C^{n+1}\to\C^2$ defined by $$F(z_1,\dots,z_{n+1})=
(z_1,z_{n+1}^3-3z_1z_{n+1}+z_2^2+\dots+z_n^2).$$ For $t=\epsilon e^{i\theta}$, the
restriction of $F$ to the hypersurface $X_t=F_1^{-1}(t)=\{t\}\times\C^n$
induces a Lefschetz fibration $F_2:X_t\to \C$, whose generic fiber
$F_2^{-1}(u)$ is the smooth hypersurface
$z_{n+1}^3-3tz_{n+1}+z_2^2+\dots+z_n^2=0$ in $\{t\}\times\C^n$. There
are two singular fibers, corresponding to the critical values
$\pm 2t^{3/2}$. It is a classical fact that the two vanishing cycles
intersect transversely in a single point ($F_2:X_t\to\C$ is a Morsification
of an $A_2$ singularity).

As $\theta$ varies from $0$ to $2\pi$, the two critical values of $F_2$ are
exchanged by the braid $\sigma^3$, where $\sigma$ is the half-twist along
the straight line segment joining the critical values.
Therefore, by definition $\rho(\sigma^3,1)$ is the monodromy of the
trivial fibration $F_1:\C^{n+1}\to\C$ around the origin, hence
$\rho(\sigma^3,1)=1$.
In other words, an $S^1$-family of symplectic Lefschetz pencils on $M$
realizing the cube of the half-twist along the arc $\gamma$ has trivial
monodromy because it bounds a $D^2$-family of $\CP^1$-valued maps on $M$,
one of which presents a degenerate critical point of type $A_2$ (this
corresponds to $t=0$ in the local model).
\end{proof}

Our last construction of elements in $\mathrm{Ker}(\rho)$ is easier to
describe in the special case of $4$-manifolds.
Assume that $\dim M=4$, and let $\delta:[0,1]\to\Sigma_0$ be an embedded
arc with end points in $N$. Let $\tau_\delta\in\mathcal{B}(\Sigma_0,N)$ be
the half-twist exchanging the two base points $\delta(0)$ and $\delta(1)$
along the arc $\delta$. Finally, let $\gamma:[0,1]\to\Delta\subset\CP^1$
be an arc with end points in $\mathcal{A}$, with a supporting pair
$(\eta',\eta'')$, and denote by $S',S''$ the corresponding
vanishing cycles.

\begin{proposition}\label{prop:fakematch}
Assume that the vanishing cycles $S',S''\subset\Sigma_0-N$ satisfy the
following properties:
$(i)$ $S'$ intersects $\delta$ in exactly one point;
$(ii)$ $S''$ is Hamiltonian isotopic to $\tau_\delta(S')$. Then
$(\sigma_\gamma,\tau_\delta)\in\Gamma(\phi)$, and $\rho(\sigma_\gamma,
\tau_\delta)=1$.
\end{proposition}

Note that, since $\tau_\delta$ represents the trivial element in the mapping
class group of $\Sigma_0$, the two vanishing cycles $S',S''$ are actually
Hamiltonian isotopic in $\Sigma_0$ (but not in $\Sigma_0-N$). In fact,
$\gamma$ is not a matching path for $\phi$, but it is a matching path for
the blown up Lefschetz fibration, and the homology class of the
associated Lagrangian sphere in the blowup $\hat{M}$ is the difference
between the exceptional classes of the blowups at the base points
$\delta(0)$ and $\delta(1)$.

\begin{proof}
The action of $\sigma_\gamma$ on $\pi_1(\Delta-\mathcal{A},\alpha_0)$ maps
$\eta'$ to $\eta''$, and $\eta''$ to $\eta''\eta'(\eta'')^{-1}$. To prove
that $(\sigma_\gamma,\tau_\delta)\in\Gamma(\phi)$, it is sufficient to check
that $L_\theta(\sigma_{\gamma *}(\eta'))=[S'']=\tau_\delta(L_\theta(\eta'))$,
which follows directly from assumption $(ii)$, and
$L_\theta(\sigma_{\gamma *}(\eta''))=
\tau_{S''}^{-1}([S'])=\tau_\delta(L_\theta(\eta''))$. In other terms, we
have to check that $\tau_{S''}^{-1}(S')$ and $\tau_\delta(S'')$ are mutually
Hamiltonian isotopic. This is easily accomplished, either by drawing a
picture, or by observing that $\tau_\delta^2$ is a Dehn twist and using the
so-called lantern relation in the mapping class group of a sphere with four
punctures to show that $\tau_{S''}(\tau_\delta^2(S'))$ is isotopic to $S'$.

A local model for the automorphism $(\sigma_\gamma,\tau_\delta)$ is
given by a family of pencils $F_t$ ($t=\epsilon e^{i\theta}$) defined over
a neighborhood of the origin in $\C^2$ by
$F_t(x,y)=(x^2+y^2-t)/x$. The pencil $F_t$ has two base points
$(0,\pm t^{1/2})$, and two critical points $(\pm i\,t^{1/2},0)$, associated
to critical values $\pm 2i\,t^{1/2}$. The smooth fibers are conics in
$\C^2$, while the singular fibers $(x\mp i\,t^{1/2})^2+y^2=0$ are unions of
two lines, each containing one of the two base points; a closer examination
shows that the roles of the two base points are exchanged when passing from
one singular fiber to the other, so that the vanishing cycles differ by a
half-twist in the fiber, as required.

When $\theta$ varies from $0$ to $2\pi$, the critical values of $F_{\epsilon
\exp(i\theta)}$ are exchanged by a half-twist; moreover, if we consider a
reference fiber $F_t^{-1}(\lambda)=\{x^2+y^2-t=\lambda x\}$ for
$|\lambda|>2\epsilon^{1/2}$, the monodromy as $\theta$ varies from $0$ to
$2\pi$ is trivial if one forgets the base points, but exchanges the two
base points $(0,\pm t^{1/2})$ by a half-twist in the fiber. Therefore,
the family of pencils $(F_{\epsilon\exp(i\theta)})$ is
indeed a local model for the situation at hand.
Since the monodromy of this family is trivial, we conclude that
$\rho(\sigma_\gamma,\tau_\delta)=1$. In other words, an $S^1$-family of
symplectic Lefschetz pencils on $M$ realizing the automorphism
$(\sigma_\gamma,\tau_\delta)$ has trivial monodromy because it bounds
a $D^2$-family of $\CP^1$-valued maps, one of which possesses a degenerate
base locus (for $t=0$, the fibers of $F_t$ intersect at the origin with
multiplicity $2$).
\end{proof}

Proposition \ref{prop:fakematch} has a natural generalization in
higher dimensions. Recall that, inside $\Sigma_0$, the base locus
$N$ represents a class Poincar\'e dual to a multiple of the
symplectic class. By varying $N$ inside a pencil, we can obtain
families of pencils on $M$ in which the base locus changes by a
Dehn twist (a transposition exchanging two points in the case
$\dim N=0$ considered above). More precisely, assume $\Sigma_0$
contains a Darboux ball $B\subset\C^{n-1}$ inside which $N$ is the
hypersurface $z_1^2+\dots+z_{n-1}^2=\epsilon$, with $0<\epsilon\ll
1$. Choose a constant $\alpha$ such that $\epsilon\ll \alpha\ll 1$
and a smooth cut-off function $\varphi:[0,\infty)\to [0,\pi]$ with
support in $[0,2\alpha]$ and equal to $\pi$ over the interval
$[0,\alpha]$. Then $(z_1,\dots,z_{n-1})\mapsto
(e^{i\varphi(|z|)}z_1,\dots, e^{i\varphi(|z|)}z_{n-1})$ is a
symplectomorphism of $\Sigma_0$, and admits a $C^1$-small
perturbation which maps $N$ to itself via a Dehn twist along the
$(n-2)$-sphere $N\cap \R^{n-1}$, thus defining an isotopy class in
$\mathcal{B}(\Sigma_0,N)$. This construction can be thought of as
a ``symplectic half-twist'' of $\Sigma_0$ along a Lagrangian disc
with boundary in $N$. As in the 4-dimensional case, there are
natural situations where the kernel of $\rho:\Gamma(\phi)\to\pi_0
\mathrm{Symp}(M,\omega)$ contains elements of the form
$(\sigma,\tau)$, where $\sigma\in B_r$ is a half-twist and $\tau$
is of the form we just described (see the next section).

\subsection{Matching paths and projective duality}

One of the main motivations for understanding ``trivial'' pencil automorphisms
(the kernel of $\rho$), besides
clarifying the relationship between pencil automorphisms and isotopy classes
of symplectomorphisms, is to optimize the search for matching paths in a
symplectic Lefschetz pencils. Indeed, we have the following obvious
statement:

\begin{proposition}\label{prop:equivrel}
If $\gamma:[0,1]\to\Delta$ is a matching path for the pencil $\phi$,
associated to a Lagrangian sphere $S_\gamma\subset M$, and if
$(b,g)\in\Gamma(\phi)$, then $b_*(\gamma)$ is also a matching path,
and the corresponding Lagrangian sphere $S_{b_*(\gamma)}$ is isotopic to
$\rho(b,g)(S_\gamma)$.
\end{proposition}

In particular, the action of the braids $b\in B_r$ for which there exists
$g\in \mathcal{B}(\Sigma_0,N)$ such that $(b,g)\in\mathrm{Ker}(\rho)$
defines an equivalence relation on the set of embedded arcs in $\Delta$
with endpoints in $\mathcal{A}$. Any arc equivalent to a given matching
path is also a matching path, and the corresponding Lagrangian spheres
are mutually isotopic.

{}From now on, we assume for simplicity that $M$ is K\"ahler. A
triple of generic holomorphic sections of a sufficiently ample
line bundle $L^{\otimes k}$, $k\gg 0$, determines a $\CP^2$-valued
map $f$ with generic local models, defined outside of a complex
codimension 3 base locus $Z\subset M$. (When $M$ is a complex
surface, $Z$ is empty and $f:M\to\CP^2$ is a branched covering).
The discriminant curve $D=\mathrm{crit}(f)\subset\CP^2$ is a
complex plane curve with cusp and node singularities. For a
generic point $p\in\CP^2$, the pencil of lines through $p$
determines a Lefschetz pencil structure on $M$ (by composition
with $f$ ); the fibers are the preimages by $f$ of the lines in
the pencil, and the singular fibers correspond to those lines
through $p$ that are tangent to the curve $D$.

Introduce the {\it dual} projective plane $(\CP^2)^*$, which is the set of
all projective lines in $\CP^2$, and let $D^*\subset(\CP^2)^*$ be the dual
curve of $D$, consisting of all the lines tangent to $D$ in $\CP^2$.
Generically the only singularities of $D^*$ are again nodes (corresponding
to lines that are tangent to $D$ in two points) and cusps (corresponding to
inflection points of $D$). A point $p^*$ in
$(\CP^2)^*$ defines a hyperplane section $\Sigma_{p^*}=
\overline{f^{-1}(L_{p^*})}\subset M$, where $L_{p^*}\subset\CP^2$ is
the line dual to $p^*$. This hyperplane section is smooth if and only if
the point $p^*$ lies outside of $D^*$. A line $\ell^*\subset(\CP^2)^*$
defines a pencil of hyperplane sections (the preimages by $f$ of the pencil
of lines through the point dual to $\ell^*$ in $\CP^2$), whose singular
fibers correspond to the points of $\ell^*\cap D^*$ (while the smooth fibers
correspond to the other points of $\ell^*$). This is a Lefschetz pencil (in
the sense of Definition \ref{def:pencil}) if and only if $\ell^*$ intersects
$D^*$ transversely at smooth points of $D^*$; otherwise the
map to $\CP^1$ corresponding to $\ell^*$ presents non-generic singularities.

If we consider a one-parameter family of lines $\ell^*_t\subset(\CP^2)^*$,
$t\in S^1$, such that each line $\ell^*_t$ is transverse to $D^*$, then we
obtain a family of Lefschetz pencils $\phi_t$, $t\in S^1$, whose monodromy
gives an element $(b,g)\in \Gamma(\phi_0)$ by considering the motion of the
critical values of $\phi_t$ as $t$ varies along $S^1$. Since the family of
maps $\phi_t$ can be thought of as a map from $M\times S^1$ to $\CP^1$
(defined outside of the base loci), it is clear that the induced symplectomorphism
of $M$ is trivial, i.e.\ $(b,g)\in\mathrm{Ker}(\rho)$.

Since the critical values of $\phi_t$ correspond to the points of
$\ell^*_t\cap D^*$, the braid $b\in B_r$ is simply the {\it braid monodromy}
of the degree $r$ plane curve $D^*$ with respect to the family of lines
$\ell^*_t$. More precisely, assume for simplicity
that the lines $\ell^*_t$ all pass through
a generic point $q_0\in (\CP^2)^*-D^*$, and consider a linear projection
$\pi:(\CP^2)^*-\{q_0\}\to\CP^1$ with pole $q_0$. Let $\Delta\subset\CP^1$
be the set of critical values of $\pi_{|D^*}$, i.e.\ the set of those fibers
of $\pi$ which pass through the singular points of $D^*$ or fail to be
transverse to $D^*$. Restricting ourselves to an affine subset
$\C\subset\CP^1$ over which the fibration $\pi$ is trivial, we can define
the braid monodromy of $D^*$,
$\psi_{D^*}:\pi_1(\C-\Delta)\to B_r$, in the following manner:
given a loop $\gamma:S^1\to\C-\Delta$, for each $t\in S^1$ the intersection
$D^*\cap \pi^{-1}(\gamma(t))$ is a configuration of $r$ points in
$\pi^{-1}(\gamma(t))\simeq
\C$; the motion of these $r$ points as $t$ varies determines a braid
$\psi_{D^*}(\gamma)\in B_r$ (see e.g.\ \cite{Teicher} or \cite{AuGo} for
more details). Now, if we consider the pencils associated to a family of lines
$\ell^*_t=\overline{\pi^{-1}(\gamma(t))}\subset(\CP^2)^*$ for some loop
$\gamma:S^1\to\C-\Delta$, then by definition we have $b=\psi_{D^*}(\gamma)$.

As a corollary of Proposition \ref{prop:equivrel} and the above
remarks, the image of any matching path for $\phi_0$ under the
action of any element of the monodromy group
$\mathrm{Im}(\psi_{D^*})\subset B_r$ is also a matching path for
$\phi_0$, and the corresponding Lagrangian spheres are mutually
isotopic. Recall that the inclusion $i:\ell^*_0-(\ell^*_0\cap D^*)
\hookrightarrow (\CP^2)^*-D^*$ induces a surjective homomorphism on
fundamental groups, and by the Zariski-van Kampen theorem, the
kernel of $i_*$ is generated by relations of the form $g\simeq
b_*(g)$ for all $g\in\pi_1(\ell^*_0-(\ell^*_0\cap D^*))$ and all
$b\in\mathrm{Im}(\psi_{D^*})\subset B_r$. The corresponding
statement for matching paths is the following: if two embedded arcs
in $\ell^*_0$ with endpoints in $\ell^*_0\cap D^*$ are isotopic as
arcs in $(\CP^2)^*$ with endpoints in $D^*$, and if one of them is a
matching path for $\phi_0$, then the other one is also a matching
path. (Another way to see this result is to consider a universal
family of hyperplane sections over $(\CP^2)^*$, and observe that
with respect to this universal fibration the notion of matching path
makes sense for arcs in $(\CP^2)^*$ with endpoints in $D^*$). Hence
the problem of classifying matching paths up to the equivalence
relation introduced at the beginning of this section reduces to the
space of isotopy classes of arcs in $(\CP^2)^*$ with endpoints in
$D^*$.

In this context, the braid monodromy of $D^*$ does not yield any
new types of elements in the kernel of
$\rho:\Gamma(\phi_0)\to\pi_0 \mathrm{Symp}(M,\omega)$, but rather
provides a geometric way of obtaining kernel elements of the form
described in Propositions \ref{prop:kerrho}--\ref{prop:fakematch}.
Namely, a node of $D^*$ corresponds to a line in $\CP^2$ that is
tangent to $D$ in two points, i.e.\ a hyperplane section of $M$
with two ordinary double points, indicative of the presence of two
mutually disjoint vanishing cycles in the pencil $\phi_0$; the
braid monodromy of $D^*$ around (the image by $\pi$ of) a node is
the square of a half-twist, and corresponds to the situation
described in Proposition \ref{prop:kerrho}~({\it a}). Similarly, a
cusp of $D^*$ corresponds to a line in $\CP^2$ that is tangent to
$D$ at an inflection point, i.e.\ a hyperplane section of $M$ with
an $A_2$ (cusp) singularity, obtained from a smooth hyperplane
section by collapsing two vanishing cycles that intersect
transversely once; the braid monodromy of $D^*$ around a cusp is
the cube of a half-twist, which corresponds to the situation
described in Proposition \ref{prop:kerrho}~({\it b}). Finally, a
line in $(\CP^2)^*$ which is tangent to $D^*$ corresponds to a
pencil of lines in $\CP^2$ passing through a point $p$ of $D$,
i.e.\ to a pencil of hyperplane sections in $M$ whose base locus
$f^{-1}(p)$ presents an ordinary double point (or, when $M$ is a
complex surface, a base point with multiplicity $2$); the braid
monodromy around a tangency of $D^*$ with the fibers of $\pi$ is a
half-twist, and corresponds to the situation discussed in
Proposition \ref{prop:fakematch}.

The above discussion should extend to the case of arbitrary
symplectic manifolds, using approximately holomorphic maps
$f:M\to\CP^2$ determined by triples of sections of $L^{\otimes k}$
(cf.\ \cite{AuGo}) satisfying suitable additional transversality
conditions. Although in this context the discriminant curve
$D\subset\CP^2$ is no longer a complex curve, a ``dual curve''
$D^*$ may still be defined by considering suitable expressions
involving the pseudoholomorphic part of the jet of the map $f$. It
follows from a general result about estimated transversality for
approximately holomorphic jets \cite{AuSe} that we can impose
conditions on the map $f$ which ensure that $D^*$ is a
well-defined symplectic curve in $(\CP^2)^*$ presenting complex
cusps and nodes of either orientation as its only singularities.
While the duality between $D$ and $D^*$ now only holds in a much
weaker sense as in the complex case, it is still reasonable to
expect that the braid monodromy of $D^*$ should give useful
information about matching paths.

To finish the discussion, we mention the following

\begin{conjecture}\label{conj:kerrho}
For pencils of sufficiently large degree $(k\gg 0)$, the kernel of
the homomorphism $\rho:\Gamma(\phi_k)\to
\pi_0\mathrm{Symp}(M,\omega)$ is generated by the three types of
elements described in Propositions
\ref{prop:kerrho}--\ref{prop:fakematch}.
\end{conjecture}

Motivation for this conjecture comes from the observation that an
$S^1$-family of {\it holomorphic} pencils on a complex projective
manifold $M\subset\CP^N$ can be described by the motion of a line
inside the dual projective space $(\CP^N)^*$, and hence is related
to the braid monodromy of the dual variety $M^*\subset(\CP^N)^*$;
however, using Lefschetz hyperplane-type arguments, one can check
that the braid monodromy of $M^*$ is generated by that of the
intersection $D^*=M^*\cap(\CP^2)^*$ for a generic plane
$(\CP^2)^*\subset(\CP^N)^*$. Also, in general, an $S^1$-family of
Lefschetz pencils on $M$ whose monodromy belongs to the identity
component in $\mathrm{Symp}(M,\omega)$ should extend to a
$D^2$-family of $\CP^1$-valued maps, in which individual members
may have singularities worse than those allowed in Lefschetz
pencils: the phenomena that are expected to occur in complex
codimension 1 are precisely those mentioned in the above
discussion.

Finally, we would like to offer the speculation that a stronger
form of Proposition \ref{prop:rhosurj} should hold for many
symplectic manifolds, e.g.\ if $\pi_0\mathrm{Symp}(M,\omega)$ is
finitely generated: namely there may exist an integer $k_0$ such
that $\rho$ is surjective for all $k\ge k_0$. In that case, by
combining Proposition \ref{prop:rhosurj} with Conjecture
\ref{conj:kerrho} one would obtain a complete description of
$\pi_0\mathrm{Symp}(M,\omega)$ in terms of pencil monodromy, and
hence reduce in principle the problem of classifying isotopy
classes of symplectomorphisms to a (probably inaccessible)
combinatorial question, similarly to what can be expected from
Theorem \ref{thm:main} for isotopy classes of Lagrangian spheres.



\begin{thebibliography}{99}

\bibitem{AMP02} {\sc J. Amoros, V. Mu\~{n}oz and F. Presas},
``Generic behavior of asymptotically holomorphic pencils'', {\it J.
Symplectic Geom.,} to appear (math.SG/0210325).

\bibitem{Au97}  {\sc D. Auroux},  ``Asymptotically holomorphic families of
symplectic submanifolds'', {\sl Geom.\ Funct.\ Anal.~\bf 7} (1997) 971--995.

\bibitem{Au00}  {\sc D. Auroux},  ``Symplectic $4$-manifolds as branched
coverings of $\CP^2$'', {\sl Invent.\ Math.~\bf 139} (2000)
551--602.

\bibitem{AuGo}  {\sc D. Auroux}, ``Symplectic maps to projective
spaces and symplectic invariants'', {\sl Turkish J.\ Math.~\bf 25} (2001)
1--42 (math.GT/0007130).

\bibitem{AuSe} {\sc D. Auroux}, ``Estimated transversality in symplectic
geometry and projective maps'', {\sl Symplectic Geometry and Mirror
Symmetry, Seoul (2000)}, World Scientific, Singapore, 2001, pp.\ 1--30
(math.SG/0010052).

\bibitem{Au01} {\sc D. Auroux, D. Gayet and J. P. Mohsen},
``Symplectic hypersurfaces in the
complement of an isotropic submanifold'', {\sl Math. Ann.~\bf 321}
(2001) 739--754.

\bibitem{AKdoubling} {\sc D. Auroux and L. Katzarkov},
``A degree doubling formula for braid monodromies and Lefschetz pencils'',
preprint.

\bibitem{Do96}  {\sc S.\,K. Donaldson}, ``Symplectic submanifolds and
almost-complex geometry'', {\sl J.~Diff.\ Geom.~\bf 44} (1996) 666--705.

\bibitem{Do00}  {\sc S.\,K. Donaldson}, ``Lefschetz pencils on symplectic
manifolds'', {\sl J.~Diff.\ Geom.~\bf 53} (1999) 205--236.

\bibitem{GM03} {\sc E. Giroux and J.\,P. Mohsen}. In preparation.

\bibitem{Gompf} {\sc R.\,E. Gompf}, ``A topological characterization of
symplectic manifolds'', math.SG/0210103.

\bibitem{MS98} {\sc D. McDuff and D. Salamon}, {\sl
Introduction to symplectic topology}, Oxford Mathematical Monographs. Clarendon
Press, Oxford University Press, New York, second edition, 1998.

\bibitem{Sei03} {\sc P. Seidel}, ``Lectures on four-dimensional Dehn
twists'', math.SG/0309012.

\bibitem{Smdoubling} {\sc I. Smith}, ``Lefschetz pencils and divisors in moduli
space'', {\sl Geom.\ Topol.}~{\bf 5} (2001), 579--608.

\bibitem{Teicher} {\sc M. Teicher}, ``Braid groups, algebraic surfaces and
fundamental groups of complements of branch curves'', {\sl Algebraic Geometry
(Santa Cruz, 1995), Proc.\ Sympos.\ Pure Math.}~{\bf 62} (part 1),
Amer.\ Math.\ Soc., Providence, 1997, pp.\ 127--150.
\end{thebibliography}
\end{document}